 \documentclass[12pt]{article}
  \usepackage{graphicx}
  \usepackage{subfig}
  \usepackage{epsfig}
  \usepackage{epstopdf}
  \usepackage{amsfonts}
  \usepackage{amsmath}
  \usepackage{mathrsfs}
  \usepackage{latexsym}
  \usepackage{color}
  \usepackage{amssymb}
  \newcommand{\ts}{\textstyle}
  \newcommand{\be}{\begin{equation}}
  \newcommand{\ee}{\end{equation}}
  \newcommand{\ba}{\begin{eqnarray}}
  \newcommand{\ea}{\end{eqnarray}}
  \newcommand{\bas}{\begin{eqnarray*}}
  \newcommand{\eas}{\end{eqnarray*}}
  \def \rla{\rangle}
  \def \lla{\langle}
  \def \beginproof{\par\noindent {\bf Proof.}\ \ }
  \def \endproof{\hskip .5cm $\Box$ \vskip .5cm}

  \def \sign{{\rm sign} }

  \def \wh{\widehat}
  \def\veps{\varepsilon}
  \def\eps{\epsilon}
  
  \def\ulamd{\underline{\lambda}}

  \def\oI{\overline{I}}

  \def\bv{\overline{v}}
  \def\bx{\overline{x}}
  
  \def\brho{\overline{\rho}}

  \def \lamd{\lambda }
  \def \ga{\gamma }
  \def \a{\alpha}

  \def\wt{\widetilde}
  
  \renewcommand{\Re}{{\rm I}\!{\rm R}}

  \newcommand{\bI}{\overline{I}}

  \newcommand{\CF}{{\cal F}}
  \newcommand{\CN}{{\cal N}}
  \newcommand{\CS}{{\cal S}}

\begin{document}

 \newtheorem{property}{Property}[section]
 \newtheorem{proposition}{Proposition}[section]
 \newtheorem{definition}{Definition}[section]
 \newtheorem{lemma}{Lemma}[section]
 \newtheorem{algorithm}{Algorithm}[section]
 \newtheorem{algorithm0}{Algorithm}
 \newtheorem{example}{Example}[section]
 \newtheorem{corollary}{Corollary}[section]
 \newtheorem{theorem}{Theorem}[section]
 \newtheorem{remark}{Remark}[section]
 \newtheorem{append}{Appendix}[section]
 \newtheorem{alemma}{Lemma}
 \newtheorem{atheorem}{Theorem}
 \newtheorem{aproposition}{Proposition}

 \begin{center}

 {\large \bf Exact penalty decomposition method for zero-norm minimization
  based on MPEC formulation
  \footnote{This work was supported by the Fundamental Research Funds for
  the Central Universities (SCUT).}}

  \vspace{1cm}

  \renewcommand{\thefootnote}{\fnsymbol{footnote}}

  Shujun Bi\footnote{Department of Mathematics,
  South China University of Technology, Guangzhou 510641,
  People's Republic of China (beamilan@163.com).},\ \
  Xiaolan Liu\footnote{Corresponding author. Department of Mathematics,
  South China University of Technology, Guangzhou 510641,
  People's Republic of China (liuxl@scut.edu.cn).}\ \ {\rm and}\ \
  Shaohua Pan\footnote{Department of Mathematics, South China University of Technology,
 Guangzhou 510641, People's Republic of China (shhpan@scut.edu.cn).}

  \vskip 1.0cm

  November 10, 2011 \\
 (First revised July 15, 2012)\\
 (Second revised March 20, 2013)\\
 (Final revision March 20, 2014)\\
  \end{center}

  \vspace{0.8cm}

  \noindent
  {\bf Abstract.}
  We reformulate the zero-norm minimization problem as an equivalent mathematical program with equilibrium constraints
  and establish that its penalty problem, induced by adding the complementarity constraint to the objective, is exact.
  Then, by the special structure of the exact penalty problem, we propose a decomposition method
  that can seek a global optimal solution of the zero-norm minimization problem under the null space condition in \cite{KXAB11}
  by solving a finite number of weighted $l_1$-norm minimization problems. To handle the weighted $l_1$-norm subproblems, 
  we develop a partial proximal point algorithm where the subproblems may be solved approximately with 
  the limited memory BFGS (L-BFGS) or the semismooth Newton-CG. Finally, we apply the exact penalty decomposition method 
  with the weighted $l_1$-norm subproblems solved by combining the L-BFGS with the semismooth Newton-CG to several types of 
  sparse optimization problems, and compare its performance with that of the penalty decomposition method \cite{ZL10}, 
  the iterative support detection method \cite{WY10} and the state-of-the-art code ${\bf FPC\_AS}$ \cite{WYGZ10}. 
  Numerical comparisons indicate that the proposed method is very efficient in terms of the recoverability 
  and the required computing time.

  \vspace{0.5cm}

  \noindent {\bf Key words:}
  zero-norm minimization; MPECs; exact penalty; decomposition method.

  \medskip

  \section{Introduction}\label{sec1}

  Let $\Re^n$ be the real vector space of dimension $n$ endowed with
  the Euclidean inner product $\lla\cdot,\cdot\rla$ and induced norm $\|\cdot\|$.
  We consider the following zero-norm minimization problem
  \ba\label{l0}
   \min_{x\in\Re^n}\Big\{\|x\|_0:\ \|Ax-b\|\le \delta\Big\},
  \ea
  where $A\in\Re^{m\times n}$ and $b\in\Re^m$ are given data,
  $\delta\ge 0$ is a given constant, and
   \[
    \|x\|_0:=\sum_{i=1}^n{\rm card}(x_i)\ \ {\rm with}\ \
    {\rm card}(x_i)=\left\{\begin{array}{ll}
                               1 & {\rm if}\ x_i\ne 0,\\
                               0 & {\rm otherwise}.
                          \end{array}\right.
  \]
  Throughout this paper, we denote by $\CF$ the feasible set of problem (\ref{l0})
  and assume that it is nonempty. This implies that (\ref{l0}) has a nonempty set
  of globally optimal solutions.

  \medskip

  The problem (\ref{l0}) has very wide applications in sparse reconstruction of signals and images
  (see, e.g., \cite{DL92,BDE09,Char07,CW08}), sparse model selection (see, e.g., \cite{Tib96,FL01}) and
  error correction \cite{CWB08}. For example, when considering the recovery of signals from noisy data,
  one may solve (\ref{l0}) with some given $\delta>0$. However, due to the discontinuity and nonconvexity of zero-norm,
  it is difficult to find a globally optimal solution of (\ref{l0}). In addition, each feasible solution of (\ref{l0})
  is locally optimal, but the number of its globally optimal solutions is finite when $\delta=0$,
  which brings in more difficulty to its solving. A common way is to obtain a favorable locally optimal solution
  by solving a convex surrogate problem such as the $l_1$-norm minimization or $l_1$-norm regularized problem.
  In the past two decades, this convex relaxation technique became very popular due to the important
  results obtained in \cite{DS89,DL92,Tib96,CDS98}. Among others, the results of \cite{DS89,DL92} quantify
  the ability of $l_1$-norm minimization problem to recover sparse reflectivity functions. For brief historical
  accounts on the use of $l_1$-norm minimization in statistics and signal processing, please see \cite{Mill02,Trop06}.
  Motivated by this, many algorithms have been proposed for the $l_1$-norm minimization or $l_1$-norm regularized problems
  (see, e.g., \cite{CR06,FNW07,KLBG07}).

  \medskip

  Observing the key difference between the $l_1$-norm and the $l_0$-norm,
  Cand\`{e}s et al. \cite{CWB08} recently proposed a reweighted $l_1$-norm
  minimization method (the idea of this method is due to Fazel \cite{FHB03} where
  she first applied it for the matrix rank minimization). This method is solving
  a sequence of convex relaxations of the following nonconvex surrogate
  \ba\label{l01}
    \min_{x\in\Re^n}\left\{\sum_{i=1}^n\ln(|x_i|+\veps):\ \|Ax-b\|\le \delta\right\}.
  \ea
  This class of surrogate problems are further studied in \cite{RSS11,ZD12}.
  In addition, noting that the $l_p$-norm $\|x\|_p^p$ tends to $\|x\|_0$ as $p\to 0$,
  many researchers seek a locally optimal solution of the problem (\ref{l0}) by solving
  the nonconvex approximation problem
  \[
    \min_{x\in\Re^n}\left\{\|x\|_p^p:\ \|Ax-b\|\le \delta\right\}
  \]
  or its regularized formulation (see, e.g., \cite{FBN07,Char07}).
  Extensive computational studies in \cite{CWB08,FBN07,Char07,RSS11}
  demonstrate that the reweighted $l_1$-norm minimization method and the $l_p$-norm
  nonconvex approximation method can find sparser solutions than the $l_1$-norm convex
  relaxation method. We see that all the methods mentioned above are developed
  by the surrogates or the approximation of zero-norm minimization problem.

  \medskip

  In this paper, we reformulate (\ref{l0}) as an equivalent MPEC (mathematical program with equilibrium constraints)
  by the variational characterization of zero-norm, and then establish that its penalty problem,
  induced by adding the complementarity constraint to the objective, is exact, i.e.,
  the set of globally optimal solutions of the penalty problem coincides with that of (\ref{l0}) when the penalty
  parameter is over some threshold. Though the exact penalty problem itself is also difficult to solve,
  we exploit its special structure to propose a decomposition method that is actually a reweighted $l_1$-norm minimization method.
  This method, consisting of a finite number of weighted $l_1$-norm minimization, is shown to yield a favorable locally
  optimal solution, and moreover a globally optimal solution of (\ref{l0}) under the null space condition in \cite{KXAB11}.
  For the weighted $l_1$-norm minimization problems, there are many softwares suitable for solving them
  such as the alternating direction method software YALL1 \cite{YALL1}, and we here propose
  a partial proximal point method where the proximal point subproblems may be solved approximately with
  the L-BFGS or the semismooth Newton-CG method (see Section 4).

  \medskip

  We test the performance of the exact penalty decomposition method with the subproblems 
  solved by combining the L-BFGS and the semismooth Newton-CG for the problems with several types of sparsity, 
  and compare its performance with that of the penalty decomposition method ({\bf QPDM}) \cite{ZL10}, 
  the iterative support detection method ({\bf ISDM}) \cite{WY10} and the state-of-the-art code
  ${\bf FPC\_AS}$ \cite{WYGZ10}. Numerical comparisons show that the proposed method has a very good robustness,
  can find the sparsest solution with desired feasibility for the Sparco collection, has comparable recoverability
  with ${\bf ISDM}$ from fewer observations for most of randomly generated problems which is higher than that of
  ${\bf FPC\_AS}$ and ${\bf QPDM}$, and requires less computing time than ${\bf ISDM}$.

  \medskip

  Notice that the {\bf ISDM} proposed in \cite{WY10} is also a reweighted $l_1$-norm minimization method
  in which, the weight vector involved in each weighted $l_1$ minimization problem is chosen as
  the support of some index set. Our exact penalty decomposition method shares this feature with {\bf ISDM},
  but the index sets to determine the weight vectors are automatically yielded by relaxing
  the exact penalty problem of the MPEC equivalent to the zero-norm problem (\ref{l0}), instead of using some heuristic strategy.
  In particular, our theoretical analysis for the exact recovery is based on the null space condition
  in \cite{KXAB11} which is weaker than the truncated null space condition in \cite{WY10}.
  Also, the weighted $l_1$-norm subproblems involved in our method are solved by combining the L-BFGS with
  the semismooth Newton-CG method, while such problems in \cite{WY10} are solved by applying YALL1 \cite{YALL1} directly.
  Numerical comparisons show that the hybrid of the L-BFGS with the semismooth Newton-CG method is effective for handling
  the weighted $l_1$-norm subproblems. The penalty decomposition method proposed by Lu and Zhang \cite{ZL10}
  aims to deal with the zero-norm minimization problem (\ref{l0}), but their method is based on a quadratic penalty
  for the equivalent augmented formulation of (\ref{l0}), and numerical comparisons show that such penalty
  decomposition method has a very worse recoverability than ${\bf FPC\_AS}$ which is designed for solving
  the $l_1$-minimization problem.

  \medskip

  Unless otherwise stated, in the sequel, we denote by $e$ a column vector
  of all $1$s whose dimension is known from the context. For any $x\in\Re^n$,
  $\sign(x)$ denotes the sign vector of $x$, $x^{\downarrow}$ is the vector
  of components of $x$ being arranged in the nonincreasing order, and $x_{I}$
  denotes the subvector of components whose indices belong to $I\subseteq\{1,\ldots,n\}$.
  For any matrix $A\in\Re^{m\times n}$, we write ${\rm Null}(A)$ as
  the null space of $A$. Given a point $x\in\Re^n$ and a constant $\ga>0$,
  we denote by $\CN(x,\ga)$ the $\ga$-open neighborhood centered at $x$.

 \section{Equivalent MPEC formulation}

  In this section, we provide the variational characterization of zero-norm
  and reformulate (\ref{l0}) as a MPEC problem. For any given $x\in\Re^n$,
  it is not hard to verify that
  \be\label{znorm1}
    \|x\|_0=\min_{v\in\Re^n}\Big\{\langle e, e-v\rangle:\
             \langle v, |x|\rangle=0,\ 0\le v\le e\Big\}.
   \ee
  This implies that the zero-norm minimization problem (\ref{l0}) is equivalent to
  \be\label{MPEC1}
    \min_{x,v\in\Re^n}\Big\{\lla e,e-v\rla:\ \|Ax-b\|\le \delta,\ \lla v,|x|\rla=0,\ 0\le v\le e\Big\},
  \ee
  which is a mathematical programming problem with the complementarity constraint:
  \[
     \lla v,|x|\rla=0,\ v\ge 0,\ |x|\ge 0.
  \]
  Notice that the minimization problem on the right hand side of (\ref{znorm1})
  has a unique optimal solution $v^*=e-\sign(|x|)$, although it is only a convex
  programming problem. Such a variational characterization of zero-norm was
  given in \cite{AA03} and \cite[Section 2.5]{Hu08}, but there it is
  not used to develop any algorithms for the zero-norm problems.
  In the next section, we will develop an algorithm for (\ref{l0})
  based on an exact penalty formulation of (\ref{MPEC1}).

  \medskip

  Though (\ref{MPEC1}) is given by expanding the original problem (\ref{l0}),
  the following proposition shows that such expansion does not increase
  the number of locally optimal solutions.
 \begin{proposition}\label{equiv-probs}
  For problems (\ref{l0}) and (\ref{MPEC1}), the following statements hold.
  \begin{description}
  \item[(a)] Each locally optimal solution of (\ref{MPEC1}) has
             the form $(x^*,e-\sign(|x^*|))$.

  \item[(b)] $x^*$ is locally optimal to (\ref{l0}) if and only if
            $(x^*,e\!-\sign(|x^*|))$ is locally optimal to (\ref{MPEC1}).
  \end{description}
  \end{proposition}
  \beginproof
  (a) Let $(x^*,v^*)$ be an arbitrary locally optimal solution of (\ref{MPEC1}).
  Then there is an open neighborhood $\CN((x^*,v^*),\ga)$ such that
  $\lla e,e-v\rla \ge \lla e,e-v^*\rla$ for all $(x,v)\in\CS\cap\CN((x^*,v^*),\ga)$,
  where $\CS$ is the feasible set of (\ref{MPEC1}). Consider (\ref{znorm1})
  associated to $x^*$, i.e.,
  \be\label{tprob1}
   \min_{v\in\Re^n}\Big\{\lla e, e-v\rla:\ \lla v, |x^*|\rla=0,\ 0\le v\le e\Big\}.
  \ee
  Let $v\in\Re^n$ be an arbitrary feasible solution of problem (\ref{tprob1})
  and satisfy $\|v-v^*\|\le \ga$. Then, it is easy to see that
  $(x^*,v)\in\CS\cap\CN((x^*,v^*),\ga)$, and so
  \(
     \lla e,e-v\rla \ge \lla e,e-v^*\rla.
  \)
  This shows that $v^*$ is a locally optimal solution of (\ref{tprob1}).
  Since (\ref{tprob1}) is a convex optimization problem,
  $v^*$ is also a globally optimal solution. However, $e-\sign(|x^*|)$ is
  the unique optimal solution of (\ref{tprob1}).
  This implies that $v^*=e-\sign(|x^*|)$.

  \medskip
  \noindent
  (b) Assume that $x^*$ is a locally optimal solution of (\ref{l0}). Then,
  there exists an open neighborhood $\CN(x^*,\ga)$ such that
  $ \|x\|_0 \ge \|x^*\|_0$ for all $x\in \CF\cap\CN(x^*,\ga)$.
  Let $(x,v)$ be an arbitrary feasible point of (\ref{MPEC1})
  such that $\|(x,v)-(x^*,e-\sign(|x^*|))\|\le \ga$. Then,
  since $v$ is a feasible point of the problem on
  the right hand side of (\ref{znorm1}), we have
 \[
   \lla e, e-v\rla \ge \|x\|_0\ge \|x^*\|_0=\lla e,e-\sign(|x^*|)\rla.
  \]
 This shows that $(x^*, e-\sign(|x^*|))$ is a locally optimal solution of (\ref{MPEC1}).

 \medskip

  Conversely, assume that $(x^*, e-\sign(|x^*|))$ is a locally optimal
  solution of (\ref{MPEC1}). Then, for any sufficiently small $\ga>0$,
  it clearly holds that $\|x\|_0\ge \|x^*\|_0$ for all $x\in\CN(x^*,\ga)$.
  This means that for any $x\in \CF\cap\CN(x^*,\ga)$, we have $\|x\|_0\ge \|x^*\|_0$.
  Hence, $x^*$ is a locally optimal solution of (\ref{l0}).
  The two sides complete the proof of part (b).
 \endproof

 \section{Exact penalty decomposition method}\label{sec3}

  In the last section we established the equivalence of (\ref{l0}) and (\ref{MPEC1}).
  In this section we show that solving (\ref{MPEC1}), and then solving (\ref{l0}),
  is equivalent to solving a single penalty problem
  \ba\label{penalty-MPEC1}
  \min_{x,v\in\Re^{n}}\Big\{\lla e,e\!-v\rla+\rho\lla v,|x|\rla:\ \|Ax-b\|\le \delta,\ 0\le v\le e \Big\}
  \ea
  where $\rho>0$ is the penalty parameter, and then develop a decomposition method
  for (\ref{MPEC1}) based on this penalty problem. It is worthwhile to point out
  that there are many papers studying exact penalty for bilevel linear programs
  or general MPECs (see, e.g., \cite{BCC89,BMS95,MPang97,LPR962}), but these references
  do not imply the exactness of penalty problem (\ref{penalty-MPEC1}).

  \medskip

  For convenience, we denote by $\CS$ and $\CS^*$ the feasible set and
  the optimal solution set of problem (\ref{MPEC1}), respectively;
  and for any given $\rho>0$, denote by $\CS_{\rho}$ and $\CS_{\rho}^*$
  the feasible set and the optimal solution set of problem (\ref{penalty-MPEC1}),
  respectively.

  \begin{lemma}\label{sol-set}
   For any given $\rho>0$, the problem (\ref{penalty-MPEC1})
   has a nonempty optimal solution set.
 \end{lemma}
 \beginproof
  Notice that the objective function of problem (\ref{penalty-MPEC1}) has
  a lower bound in $\CS_\rho$, to say $\a^*$. Therefore, there must exist
  a sequence $\{(x^k,v^k)\}\subset\CS_{\rho}$ such that for each $k$,
  \be\label{temp-equa1}
    \lla e,e-v^k\rla +\rho\lla |x^k|,v^k\rla \le \a^*+\frac{1}{k}.
  \ee
  Since the sequence $\{v^k\}$ is bounded, if the sequence $\{x^k\}$
  is also bounded, then by letting $(\bx,\bv)$ be an arbitrary limit point
  of $\{(x^k,v^k)\}$, we have $(\bx,\bv)\in\CS_\rho$ and
  \[
    \lla e,e-\bv\rla + \rho\lla |\bx|,\bv\rla\le \a^*,
  \]
  which implies that $(\bx,\bv)$ is a globally optimal solution of (\ref{penalty-MPEC1}).
  We next consider the case where the sequence $\{x^k\}$ is unbounded. Define
  the disjoint index sets $I$ and $\bI$ by
  \[
    I:=\left\{i\in\{1,\ldots,n\}\ |\ \{x_i^k\}\ {\rm is\ unbounded}\right\}
    \ \ {\rm and}\ \
    \bI:=\{1,\ldots,n\}\backslash I.
  \]
  Since $\{v^k\}$ is bounded, we without loss of generality assume that
  it converges to $\bv$. From equation (\ref{temp-equa1}), it then follows that
  $\bv_I=0$. Note that the sequence $\{x^k\}\subset\Re^{n}$ satisfies
 $A_{I}x_{I}^k+A_{\bI}x_{\bI}^k=b+\Delta^k$ with $\|\Delta^k\|\le \delta$.
 Since the sequences $\{x_{\bI}^k\}$ and $\{\Delta^k\}$ are bounded,
 we may assume that they converge to $\bx_{\bI}$ and $\Delta$, respectively.
 Then, from the closedness of the set $A_{I}\Re^I$, there exists an $\xi\in\Re^{I}$
 such that $A_{I}\xi+A_{\bI}\bx_{\bI}=b+\Delta$ with $\|\Delta\|\le \delta$.
 Letting $\bx=(\xi,\bx_{\bI})\in\Re^n$, we have $(\bx,\bv)\in \CS_\rho$.
 Moreover, the following inequalities hold:
 \bas
  \lla e,e-\bv\rla+\rho\lla |\bx|,\bv\rla
  = \lla e,e-\bv\rla + \rho\sum_{i\in \bI} |\bx_i|\bv_i
  = \lim_{k\rightarrow\infty}\Big[\lla e,e-v^k\rla + \rho\sum_{i\in \bI}\lla |x_i^k|,v_i^k\rla\Big]\\
  \le  \lim_{k\rightarrow\infty}\Big[\lla e,e-v^k\rla +  \rho\sum_{i\in I}\lla |x_i^k|,v_i^k\rla
      + \rho\sum_{i\in \bI}\lla |x_i^k|,v_i^k\rla\Big]\le\a^*,\quad
 \eas
 where the first equality is using $\bv_I=0$, and the last equality
 is due to (\ref{temp-equa1}). This shows that $(\bx,\bv)$ is a global
 optimal solution of (\ref{penalty-MPEC1}). Thus, we prove that for any given $\rho>0$,
 the problem (\ref{penalty-MPEC1}) has a nonempty set of globally optimal solutions.
 \endproof

  To show that the solution of problem (\ref{MPEC1}) is equivalent
  to that of a single penalty problem (\ref{penalty-MPEC1}), i.e.,
  to prove that there exists $\brho>0$ such that the set of global
  optimal solutions of (\ref{MPEC1}) coincides with that of
  (\ref{penalty-MPEC1}) with $\rho>\brho$, we need the following lemma.
  Since this lemma can be easily proved by contradiction, we here do not
  present its proof.

 \begin{lemma}\label{lemma-l00}
  Given $M\in\Re^{m\times n}$ and $q\in\Re^{m}$.
  If
  \(
     r=\min\big\{\|z\|_0: \|Mz-q\|\le \delta\big\}>0,
  \)
  then there exists $\a>0$ such that for all $z$ with $\|Mz-q\|\leq \delta$,
  we have $|z|_{r}^{\downarrow}>\a$, where $|z|_r^{\downarrow}$ means
  the $r$th component of $|z|^{\downarrow}$ which is the vector of components
  of $|z|\in\Re^n$ being arranged in the nonincreasing order
  $|z|_1^{\downarrow}\ge |z|_2^{\downarrow}\ge\cdots\ge|z|_n^{\downarrow}$.
 \end{lemma}

  Now we are in a position to establish that (\ref{penalty-MPEC1})
  is an exact penalty problem of (\ref{MPEC1}).

 \begin{theorem}\label{exact-penalty}
  There exists a constant $\brho\!>0$ such that $\CS^*$
  coincides with $\CS_{\rho}^*$ for all $\rho\!>\brho$.
 \end{theorem}
 \beginproof
  Let $r=\min\{\|x\|_0: \|Ax-b\|\le \delta\}$. We only need to consider
  the case where $r>0$ (if $r=0$, the conclusion clearly holds for all $\rho>0$).
  By Lemma \ref{lemma-l00}, there exists $\a>0$ such that for all $x$
  satisfying $\|Ax-b\|\le \delta$, it holds that $|x|_r^\downarrow>\a$.
  This in turn means that $|x|_r^\downarrow>\a$ for all $(x,v)\in\CS$
  and $(x,v)\in\CS_{\rho}$ with any $\rho>0$.

  \medskip

  Let $(\bx,\bv)$ be an arbitrary point in $\CS^*$. We prove that
  $(\bx,\bv)\in \CS_{\rho}^*$ for all $\rho>1/\a$, and then $\brho=1/\a$
  is the one that we need. Let $\rho$ be an arbitrary constant with
  $\rho>1/\a$. Since $(\bx,\bv)\in \CS^*$, from Proposition \ref{equiv-probs}(a)
  and the equivalence between (\ref{MPEC1}) and (\ref{l0}), it follows that
  $\bv=e-\sign(|\bx|)$ and $\|\bx\|_0=r$.
  Note that $|x|_r^\downarrow>\a>1/\rho$ for any $(x,v)\in \CS_{\rho}$
  since $\rho>1/\alpha$. Hence, for any $(x,v)\in \CS_{\rho}$,
  the following inequalities hold:
 \bas
   \lla e,e-v\rla +\rho\lla v,|x|\rla
   &=& \sum_{i=1}^{n}(1-v_i+\rho v_i|x|_i)
    \geq \sum_{|x|_i> 1/\rho}(1-v_i+\rho v_i|x|_i)\\
   &\ge& r =\lla e,e-\bv\rla + \rho\lla \bv,|\bx|\rla
 \eas
  where the first inequality is using $1-v_i+\rho v_i|x|_i\geq0$
  for all $i\in\{1,2,\ldots,n\}$, and the second inequality is since
  $|x|_r^\downarrow>1/\rho$ and $0\leq v_i\leq 1$. Since $(x,v)$
  is an arbitrary point in $\CS_{\rho}$ and $(\bx,\bv)\in\CS_{\rho}$,
  the last inequality implies that $(\bx,\bv)\in \CS_{\rho}^*$.

 \medskip

  Next we show that if $(\bx,\bv)\in \CS_{\rho}^*$ for $\rho>1/\alpha$,
  then $(\bx,\bv)\in \CS^*$. For convenience, let
 \[
  I_{-}:=\left\{i\in\{1,\ldots,n\}\ |\ |\bx|_i\le \rho^{-1}\right\}
  \ \ {\rm and}\ \
  I_{+}:=\left\{i\in\{1,\ldots,n\}\ |\ |\bx|_i> \rho^{-1}\right\}.
 \]
 Note that $1-\bv_i+\rho \bv_i|\bx|_i\ge 0$ for all $i$,
 and $1-\bv_i+\rho \bv_i|\bx|_i\ge 1$ for all $i\in I_{+}$.
 Also, the latter together with $|\bx|_r^{\downarrow}>\a>1/\rho$
 implies $|I_{+}|\ge r$. Then, for any $(x,v)\in\CS$ we have
  \bas
    \lla e,e-v\rla&=& \lla e,e-v\rla + \rho\lla v,|x|\rla
    \ge\lla e,e-\bv\rla + \rho\lla \bv,|\bx|\rla\\
    &=&\sum_{i=1}^{n}(1-\bv_i+ \rho\bv_i|\bx|_i)
    \geq \sum_{i\in I_{+}}(1-\bv_i+\rho \bv_i|\bx|_i)\geq r,
   \eas
 where the first inequality is using $\CS \subseteq \CS_{\rho}$.
 Let $\wt{x}\in\Re^{n}$ be such that $\|\wt{x}\|_0=r$
 and $\|A\wt{x}-b\|\le\delta$. Such $\wt{x}$ exists by the definition
 of $r$. Then $(\wt{x},e-\sign(|\wt{x}|))\in \CS\subseteq \CS_{\rho}$, and
 from the last inequality
 \(
  r=\lla e,e-(e-\sign(|\wt{x}|))\rla\ge\lla e,e-\bv\rla + \rho\lla \bv,|\bx|\rla\ge r.
 \)
 Thus,
  \[
   r=\lla e,e-\bv\rla + \rho\lla\bv,|\bx|\rla
   =\sum_{i\in I_{-}}(1-\bv_i+\rho\bv_i|\bx|_i)
    +\sum_{i\in I_{+}}(1-\bv_i+\rho\bv_i|\bx|_i).
  \]
 Since $1-\bv_i+\rho\bv_i|\bx_i|\ge 0$ for $i\in I_{-}$,
 $1-\bv_i+\rho\bv_i|\bx|_i \ge 1$ for $i\in I_{+}$ and $|I_{+}|\ge r$,
 from the last equation it is not difficult to deduce that
  \be\label{tequa2}
  1-\bv_i+\rho \bv_i|\bx|_i=0\ \ {\rm for}\ \ i\in I_{-},\ \
  1-\bv_i+\rho \bv_i|\bx|_i=1\ \ {\rm for}\ \ i\in I_{+},\ \ {\rm and}\ \
  |I_{+}|=r.
 \ee
 Since each  $1-\bv_i+\rho \bv_i|\bx|_i$ is nonnegative,
 the first equality in (\ref{tequa2}) implies that
 \[
    \bv_i=1\ \ {\rm and}\ \ |\bx|_i=0\ \ {\rm for}\ \ i\in I_{-};
 \]
 while the second equality in (\ref{tequa2}) implies that
 \(
    \bv_i=0
 \)
 for $i\in I_{+}$. Combining the last equation with $|I_{+}|=r$,
 we readily obtain $\|\bx\|_0=r$ and $\bv=e-\sign(|\bx|)$.
 This, together with $\|A\bx-b\|\le\delta$, shows that $(\bx,\bv)\in \CS^*$.
 Thus, we complete the proof.
 \endproof

 Theorem \ref{exact-penalty} shows that to solve the MPEC problem (\ref{MPEC1}),
 it suffices to solve (\ref{penalty-MPEC1}) with a suitable large $\rho>\brho$.
 Since the threshold $\brho$ is unknown in advance, we need to solve
 a finite number of penalty problems with increasing $\rho$.
 Although problem (\ref{penalty-MPEC1}) for a given $\rho>0$
 has a nonconvex objective function, its separable structure leads to
 an explicit solution with respect to the variable $v$ if the variable
 $x$ is fixed. This motivates us to propose the following decomposition method
 for (\ref{MPEC1}), and consequently for (\ref{l0}).

 \begin{algorithm} \label{Alg1}(Exact penalty decomposition method for problem (\ref{MPEC1}))
 \begin{description}
 \item[(S.0)] Given a tolerance $\epsilon>0$ and a ratio $\sigma>1$. Choose
              an initial penalty parameter \hspace*{0.1cm} $\rho_0>0$ and
              a starting point $v^0=e$. Set $k:=0$.

 \item[(S.1)] Solve the following weighted $l_1$-norm minimization problem
              \be\label{subprob11}
                 x^{k+1}\in \mathop{\arg\min}_{x\in\Re^{n}}\Big\{\lla v^k,|x|\rla:\ \|Ax-b\|\le \delta\Big\}.
              \ee

 \item[(S.2)] If $|x_i^{k+1}|> 1/{\rho_k}$, then set $v_i^{k+1}=0$; and otherwise
              set $v_i^{k+1}=1$.

  \item[(S.3)] If $\lla v^{k+1},|x^{k+1}|\rla\le\epsilon$, then stop;
               and otherwise go to (S.4).

 \item[(S.4)] Let $\rho_{k+1}:=\sigma\rho_k$ and $k:=k+1$, and then go to Step (S.1).

 \end{description}
 \end{algorithm}

 \begin{remark}\label{remark-alg1}
   Note that the vector $v^{k+1}$ in (S.2) is an optimal solution of the problem
  \be\label{subprob1}
      \min_{0\le v\le e}\Big\{\lla e,e-v\rla + \rho_k\lla |x^{k+1}|,v\rla\Big\}.
  \ee
  So, Algorithm \ref{Alg1} is solving the nonconvex penalty
  problem (\ref{penalty-MPEC1}) in an alternating way.
%
 \end{remark}

  The following lemma shows that the weighted $l_1$-norm subproblem (\ref{subprob11})
  has a solution.

 \begin{lemma}\label{sol-subprob1}
  For each fixed $k$, the subproblem (\ref{subprob11}) has an optimal solution.
 \end{lemma}
 \beginproof
  Note that the feasible set of (\ref{subprob11}) is $\CF$, which is nonempty
  by the given assumption in the introduction, and its objective function is bounded below
  in the feasible set. Let $\nu^*$ be the infinum of the objective function
  of (\ref{subprob11}) on the feasible set. Then there exists a feasible sequence
  $\{x^l\}$ such that $\lla v^k,|x^l|\rla \to \nu^*$ as $l\to\infty$.
  Let $I:=\{i\ |\ v_i^k=1\}$. Then, noting that
  $\lla v^k,|x^l|\rla=\sum_{i\in I}|x_{i}^{l}|$, we have that the sequence
  $\{x_{I}^{l}\}$ is bounded. Without loss of generality, we assume that
  $x_{I}^{l} \to \wt{x}_{I}$. Let $y^l=Ax^l -b$. Noting that $\|y^l\|\le \delta$,
  we may assume that $\{y^l\}$ converges to $\wt{y}$. Since the set
  $A_{\oI}\Re^{|\oI|}$ is closed and $A_{\oI}x_{\oI}^l = y^l-A_{I}x_{I}^l+b$
  for each $l$, where $\oI=\{1,2,\ldots,n\}\backslash I$, there exists
  $\wt{x}_{\oI}\in\Re^{|\oI|}$ such that
  $A_{\oI}\wt{x}_{\oI} = \wt{y}-A_{I}\wt{x}_{I}+b$, i.e.,
  $A_{\oI}\wt{x}_{\oI}+A_{I}\wt{x}_{I}-b=\wt{y}$.
  Let $\wt{x}=(\wt{x}_{I};\wt{x}_{\oI})$. Then, $\wt{x}$ is
  a feasible solution to (\ref{subprob11}) with
  \(
     \lla v^k,\wt{x}\rla=\nu^*.
  \)
  So, $\wt{x}$ is an optimal solution of (\ref{subprob11}).
 \endproof

  For Algorithm \ref{Alg1}, we can establish the following finite termination result.

 \begin{theorem}\label{fstop-Alg1}
  Algorithm \ref{Alg1} will terminate after at most
  $\lceil\frac{\ln (n)-\ln(\epsilon\rho_0)}{\ln \sigma}\rceil$ iterations.
 \end{theorem}
 \beginproof
  By Lemma \ref{sol-subprob1}, for each $k\ge 0$ the subproblem (\ref{subprob11})
  has a solution $x^{k+1}$. From Step (S.2) of Algorithm \ref{Alg1},
  we know that $v_{i}^{k+1}=1$ for those $i$ with $|x_i^{k+1}|\le 1/\rho_k$.
  Then,
  \[
    \ts\lla v^{k+1},|x^{k+1}|\rla=\sum_{\{i:\ v_i^{k+1}=1\}}|x_i^{k+1}|
    \le \frac{n}{\rho_k}.
 \]
 This means that, when $\rho_k\ge \frac{n}{\epsilon}$,
 Algorithm \ref{Alg1} must terminate. Note that $\rho_k\ge \sigma^k\rho_0$.
 Therefore, Algorithm \ref{Alg1} will terminate when
 $\sigma^k\rho_0\ge \frac{n}{\epsilon}$, i.e.,
 $k\ge \lceil\frac{\ln (n)-\ln(\epsilon\rho_0)}{\ln \sigma}\rceil$.
 \endproof

  We next focus on the theoretical results of Algorithm \ref{Alg1} for the case
  where $\delta=0$. To this end, let $x^*$ be an optimal solution of the zero-norm
  problem (\ref{l0}) and write
  \[
    I^*\!=\big\{i\ |\ x_i^* \neq 0\big\}\ \ {\rm and}\ \
    \oI^*=\!\{1,\ldots,n\}\backslash\!I^*.
  \]
  In addition, we also need the following null space condition for a given vector $v\in\Re_{+}^n$:
  \be\label{PNSC}
   \lla v_{I^*},|y_{I^*}|\rla<\lla v_{\oI^*},|y_{\oI^*}|\rla\ \ {\rm for\ any}\
   0\neq y\in {\rm Null}(A).
 \ee

 \begin{theorem}\label{NSC-Alg1}
  Assume that $\delta=0$ and $v^k$ satisfies the condition (\ref{PNSC}) for some
  nonnegative integer $k$.
  Then, $x^{k+1}=x^*$. If, in addition, $v^{k+1}$ also satisfies
  the condition (\ref{PNSC}), then the vector $v^{k+l}$ for all $l\ge 2$ satisfy
  the null space condition (\ref{PNSC}) and $x^{k+l+1}=x^*$ for all $l\ge 1$.
  Consequently, if $v^0$ satisfies the condition (\ref{PNSC}), then $x^{k}=x^*$ for all $k\ge 1$.
 \end{theorem}
 \beginproof
  We first prove the first part. Suppose that $x^{k+1}\ne x^*$. Let $y^{k+1}=x^{k+1}-x^*$.
  Clearly, $0\ne y^{k+1}\in{\rm Null}(A)$. Since $v^k$ satisfies the condition (\ref{PNSC}),
  we have that
  \begin{align}\label{temp-ineq1}
    \lla v_{\overline{I}^*}^k, |y_{\overline{I}^*}^{k+1}|\rla
   >\lla v_{I^*}^k,|y_{I^*}^{k+1}|\rla.
  \end{align}
  On the other hand, from step (S.1), it follows that
  \[
   \lla v^k,|x^*|\rla\ge \lla v^k,|x^{k+1}|\rla
   =\lla v^k,|x^*+y^{k+1}|\rla \ge \lla v^k,|x^*|\rla
    -\lla v_{I^*}^k,|y_{I^*}^{k+1}|\rla +\lla v_{\overline{I}^*}^k,|y_{\overline{I}^*}^{k+1}|\rla.
 \]
  This implies that $\lla v_{I^*}^k,|y_{I^*}^{k+1}|\rla\ge\lla v_{\overline{I}^*}^k,|y_{\overline{I}^*}^{k+1}|\rla$.
  Thus, we obtain a contradiction to (\ref{temp-ineq1}). Consequently, $x^{k+1}=x^*$.
  Since $v^{k+1}$ also satisfies the null space condition (\ref{PNSC}),
  using the same arguments yields that $x^{k+2}=x^*$.
  We next show by induction that $v^{k+l}$ for all $l\ge 2$ satisfy the condition (\ref{PNSC})
  and $x^{k+l+1}=x^*$ for all $l\ge 2$. To this end, we define
  \[
    I_{\nu}:=\left\{i\ |\ x_i^* >\frac{1}{\rho_{k+\nu}}\right\}\ \ {\rm and}\ \
    \overline{I}_{\nu}:=\left\{i\ |\ 0<x_i^* \leq\frac{1}{\rho_{k+\nu}}\right\}
  \]
  for any given nonnegative integer $\nu$. Clearly, $I^*=I_{\nu}\cup\overline{I}_\nu$.
  Also, by noting that $\rho_{k+\nu+1}=\sigma\rho_{k+\nu}$ by (S.4) and $\sigma>1$,
  we have that $I_{\nu}\subseteq I_{\nu+1}\subseteq I^*$.
  We first show that the result holds for $l=2$. Since $x^{k+1}=x^*$,
  we have $v_{I_0}^{k+1}=0$ and $v^{k+1}_{\overline{I}_0}=e$ by step (S.2).
  Since $x^{k+2}=x^*$, from step (S.2) it follows that $v_{I_1}^{k+2}=0$
  and $v^{k+2}_{\overline{I}_1}=e$. Now we obtain that
  \begin{align}\label{ineq-obj1}
   \langle v^{k+2}_{I^{*}}, |y_{I^*}|\rangle
   &=\langle v^{k+2}_{I_1}, |y_{I_1}|\rangle + \langle v^{k+2}_{\overline{I}_1}, |y_{\overline{I}_1}|\rangle
   =\langle v^{k+2}_{\overline{I}_1}, |y_{\overline{I}_1}|\rangle
    \le \langle v^{k+1}_{\overline{I}_0}, |y_{\overline{I}_0}|\rangle\nonumber\\
   &=\langle v^{k+1}_{I_0}, |y_{I_0}|\rangle + \langle v^{k+1}_{\overline{I}_0}, |y_{\overline{I}_0}|\rangle
    =\langle v^{k+1}_{I^*},|y_{I^*}|\rangle\nonumber\\
   &<\lla v^{k+1}_{\oI^*},|y_{\oI^*}|\rla =\langle v^{k+2}_{\overline{I}^{*}},|y_{\oI^*}|\rla
  \end{align}
  for any $0\ne y\in{\rm Null}(A)$, where the first equality is due to
  $I^*=I_{1}\cup\overline{I}_1$, the second equality is using $v_{I_1}^{k+2}=0$,
  the first inequality is due to $\overline{I}_1\subseteq \overline{I}_0$,
  $v^{k+2}_{\overline{I}_1}=e$ and $v^{k+1}_{\overline{I}_0}=e$, the second inequality
  is using the assumption that $v^{k+1}$ satisfies the null space condition (\ref{PNSC}),
  and the last equality is due to $v^{k+1}_{\overline{I}^*}=e$ and $v^{k+2}_{\overline{I}^*}=e$.
  The inequality (\ref{ineq-obj1}) shows that $v^{k+2}$ satisfies the null space condition
  (\ref{PNSC}), and using the same arguments as for the first part yields that $x^{k+3}=x^*$.
  Now assuming that the result holds for $l(\ge 2)$, we show that it holds for $l+1$.
  Indeed, using the same arguments as above, we obtain that
  \begin{align}\label{ineq-obj}
   \langle v^{k+l+1}_{I^{*}}, |y_{I^*}|\rangle
   &=\langle v^{k+l+1}_{I_l}, |y_{I_l}|\rangle + \langle v^{k+l+1}_{\overline{I}_l}, |y_{\overline{I}_l}|\rangle
   =\langle v^{k+l+1}_{\overline{I_l}}, |y_{\overline{I}_l}|\rangle
    \le \langle v^{k+l}_{\overline{I}_{l-1}}, |y_{\overline{I}_{l-1}}|\rangle\nonumber\\
   &=\langle v^{k+l}_{I_{l-1}}, |y_{I_{l-1}}|\rangle + \langle v^{k+l}_{\overline{I}_{l-1}}, |y_{\overline{I}_{l-1}}|\rangle
    =\langle v^{k+l}_{I^*},|y_{I^*}|\rangle\nonumber\\
   &<\lla v^{k+l}_{\oI^*},|y_{\oI^*}|\rla =\langle v^{k+l+1}_{\overline{I}^{*}},|y_{\oI^*}|\rla
  \end{align}
  for any $0\ne y\in{\rm Null}(A)$. This shows that $v^{k+l+1}$ satisfies the null space condition
  (\ref{PNSC}), and using the same arguments as the first part yields that $x^{k+l+2}=x^*$.
  Thus, we show that $v^{k+l}$ for all $l\ge 2$ satisfy the condition (\ref{PNSC})
  and $x^{k+l+1}=x^*$ for all $l\ge 2$.
 \endproof

  Theorem \ref{NSC-Alg1} shows that, when $\delta=0$, if there are two successive
  vectors $v^{k}$ and $v^{k+1}$ satisfy the null space condition (\ref{PNSC}), then
  the iterates after $x^{k}$ are all equal to some optimal solution of (\ref{l0}).
  Together with Theorem \ref{fstop-Alg1}, this means that Algorithm \ref{Alg1}
  can find an optimal solution of (\ref{l0}) within a finite number of iterations
  under (\ref{PNSC}). To the best of our knowledge, the condition (\ref{PNSC}) is first
  proposed by Khajehnejad et al. \cite{KXAB11}, which generalizes the null space condition of \cite{SXY08}
  to the case of weighted $l_1$-norm minimization and is weaker than the truncated null space condition
  \cite[Defintion 1]{WY10}.

%

 \section{Solution of weighted $l_1$-norm subproblems}

  This section is devoted to the solution of the subproblems involved
  in Algorithm \ref{Alg1}:
  \be\label{weight-l1}
   \min_{x\in\Re^n}\left\{\lla v,|x|\rla:\ \|Ax-b\|\le \delta\right\}
  \ee
  where $v\in\Re^n$ is a given nonnegative vector. For this problem, one may reformulate it as
  a linear programming problem (for $\delta=0$) or a second-order cone programming problem (for $\delta>0$),
  and then directly apply the interior point method software SeDuMi \cite{Stur01} or $l_1$-MAGIC \cite{CR06} for solving it.
  However, such second-order type methods are time-consuming, and are not suitable for handling large-scale problems.
  Motivated by the recent work \cite{JST12}, we in this section develop a partial proximal point algorithm (PPA) for
  the following reformulation of the weighted $l_1$-norm problem (\ref{weight-l1}):
  \be\label{weight1-l2}
     \min_{u\in\Re^m,x\in\Re^n}\left\{\lla v,|x|\rla + \frac{\beta}{2}\|u\|^2: Ax+u=b\right\}
     \ \ {\rm for\ some}\ \beta>0.
  \ee
  Clearly, (\ref{weight1-l2}) is equivalent to (\ref{weight-l1}) if $\delta>0$; and otherwise
  is a penalty problem of (\ref{weight-l1}).

  \medskip

  Given a starting point $(u^0,x^0)\in\Re^m\times\Re^n$, the partial PPA for (\ref{weight1-l2})
  consists of solving approximately a sequence of strongly convex minimization problems
  \be\label{subpro-PPA}
   (u^{k+1},x^{k+1}) \approx
   \arg\min_{u\in\Re^m,x\in\Re^n}\left\{\lla v,|x|\rla + \frac{\beta}{2}\|u\|^2+\frac{1}{2\lamd_k}\big\|x-x^k\big\|^2: Ax+u=b\right\},
  \ee
  where $\{\lambda_k\}$ is a sequence of parameters satisfying
  $0<\lambda_k\uparrow \overline{\lambda}\le +\infty$. For the global and
  local convergence of this method, the interested readers may refer to
  Ha's work \cite{Ha90}, where he first considered such PPA for finding
  a solution of generalized equations. Here we focus on the approximate
  solution of the subproblems (\ref{subpro-PPA}) via the dual method.

  \medskip

  Let $L:\Re^m\times\Re^n\times\Re^m\to\Re$ denote the Lagrangian function
  of the problem (\ref{subpro-PPA})
  \[
    L(u,x,y):= \lla v,|x|\rla + \frac{\beta}{2}\|u\|^2+\frac{1}{2\lamd_k}\|x-x^k\|^2
               +\lla y,Ax+u-b\rla.
  \]
  Then the minimization problem (\ref{subpro-PPA}) is expressed as
  \(
    \min_{(u,x)\in\Re^{m}\times\Re^n}\sup_{y\in\Re^m}L(u,x,y).
  \)
  Also, by \cite[Corollary 37.3.2]{Roc70} and the coercivity of $L$
  with respect to $u$ and $x$,
  \be\label{dual-gap}
    \min_{(u,x)\in\Re^{m}\times\Re^n}\sup_{y\in\Re^m}L(u,x,y)
    =\sup_{y\in\Re^m}\min_{(u,x)\in\Re^{m}\times\Re^n}L(u,x,y).
  \ee
  This means that there is no dual gap between the problem (\ref{subpro-PPA}) and its dual problem
  \be\label{dual-PPA}
   \sup_{y\in\Re^m}\min_{(u,x)\in\Re^{m}\times\Re^n}L(u,x,y).
  \ee
  Hence, we can obtain the approximate optimal solution $(u^{k+1},x^{k+1})$
  of (\ref{subpro-PPA}) by solving (\ref{dual-PPA}).
  To give the expression of the objective function of (\ref{dual-PPA}),
  we need the following operator
  $S_{\lamd}(\cdot,v)\!:\Re^n\to\Re^n$ associated to the vector $v$ and any $\lamd>0$:
  \[
    S_{\lamd}(z,v):=\arg\min_{x\in\Re^n}\left\{\lla v,|x|\rla + \frac{1}{2\lamd}\|x-z\|^2\right\}.
  \]
  An elementary computation yields the explicit expression of the operator $S_{\lamd}(\cdot,v)$:
  \be\label{Slamd}
    S_{\lamd}(z,v) ={\rm sign}(z)\odot{\rm max}\big\{|z|-\lamd v,0\big\}
    \quad\ \forall z\in\Re^n,
  \ee
  where ``$\odot$" means the componentwise product of two vectors,
  and for any $z\in\Re^n$,
  \ba\label{term-equa}
   \lla v,|S_{\lamd}(z,v)|\rla =\left\lla{\rm sign}(z)\odot v,S_{\lamd}(z,v)\right\rla,\qquad\qquad\nonumber\\
    \lla S_{\lamd}(z,v),S_{\lamd}(z,v)\rla
    =\lla S_{\lamd}(z,v),z\rla-\lamd\left\lla{\rm sign}(z)\odot v,S_{\lamd}(z,v)\right\rla.
  \ea
  From the definition of $S_{\lamd}(\cdot,v)$ and equation (\ref{term-equa}),
  we immediately obtain that
  \bas
   \min_{(u,x)\in\Re^{m}\times\Re^n}L(u,x,y)
   = -b^Ty -\frac{1}{2\beta}\|y\|^2-\frac{1}{2\lamd_k}\big\|S_{\lamd_k}(x^k-\lamd_k A^Ty,v)\big\|^2
     + \frac{1}{2\lamd_k}\|x^k\|^2.
  \eas
  Consequently, the dual problem (\ref{dual-PPA}) is equivalent to
  the following minimization problem
  \be\label{dual-PPA1}
    \min_{y\in\Re^m} \Phi(y):=b^Ty +\frac{1}{2\beta}\|y\|^2+\frac{1}{2\lamd_k}\big\|S_{\lamd_k}(x^k-\lamd_k A^Ty,v)\big\|^2.
  \ee

  The following lemma summarizes the favorable properties of the function $\Phi$.

  \begin{lemma}\label{pro-thetak}
   The function $\Phi$ defined by (\ref{dual-PPA1}) has the following properties:
   \begin{description}
   \item[(a)] $\Phi$ is a continuously differentiable convex function with
              gradient given by
              \[
                \nabla \Phi(y) = b+\beta^{-1}y - AS_{\lamd_k}(x^k-\lamd_k A^Ty,v)
                \quad\ \forall y\in\Re^m.
              \]

   \item[(b)] If $\wh{y}^k$ is a root to the system $\nabla \Phi(y)=0$,
              then $(\wh{u}^{k+1},\wh{x}^{k+1})$ defined by
              \[
                 \wh{u}^{k+1}:=-\beta^{-1}\wh{y}^{k}\ \ {\rm and}\ \
                 \wh{x}^{k+1}:=S_{\lamd_k}(x^k-\lamd_k A^T\wh{y}^k,v)
              \]
              is the unique optimal solution of the primal problem (\ref{subpro-PPA}).

   \item[(c)] The gradient mapping $\nabla\Phi(\cdot)$ is Lipschitz continuous
              and strongly semismooth.

   \item[(d)] The Clarke's generalized Jacobian of the mapping $\nabla\Phi$
              at any point $y$ satisfies
              \ba\label{upper-estimate}
                 \partial(\nabla\Phi)(y)
                 &\subseteq& \beta^{-1} I +\lamd_k A\partial_x H(z,v)A^T
                 :=\wh{\partial}^2\Phi(y).
              \ea
              where $z=x^k-\lamd_k A^Ty$ and $H(x,v):={\rm sign}(x)\odot\max\{|x|-\lamd_kv,0\}$.
   \end{description}
  \end{lemma}
  \beginproof
  (a) By the definition of $S_{\lamd}(\cdot,v)$ and equation (\ref{term-equa}),
      it is not hard to verify that
  \bas
    \frac{1}{2\lamd}\|S_{\lamd}(z,v)\|^2
    =\frac{1}{2\lamd}\|z\|^2 - \min_{x\in\Re^n}\left\{\lla v,|x|\rla
      + \frac{1}{2\lamd}\|x-z\|^2\right\}.
  \eas
  Note that the second term on the right hand side is the Moreau-Yosida
  regularization of the convex function $f(x):=\lla v,|x|\rla$. From \cite{Roc70}
  it follows that $\|S_{\lamd}(\cdot,v)\|^2$ is continuously differentiable,
  which implies that $\Phi$ is continuously differentiable.

  \medskip
  \noindent
  (b) Note that $\widehat{y}^k$ is an optimal solution of (\ref{dual-PPA})
  and there is no dual gap between the primal problem (\ref{subpro-PPA})
  and its dual (\ref{dual-PPA}) by equation (\ref{dual-gap}).
  The desired result then follows.

  \medskip
  \noindent
  (c) The result is immediate by the expression of $\Phi$ and
  $S_{\lamd_k}(\cdot,v)$.

  \medskip
  \noindent
  (d) The result is implied by the corollary in \cite[p.75]{Clar83}.
  Notice that the inclusion in (\ref{upper-estimate}) can not be
  replaced by the equality since $A$ is assumed to be of full row rank.
  \endproof
  
  \begin{remark}\label{Hsubdiff}
   For a given $v\in\mathbb{R}_{+}^n$, from \cite[Chaper 2]{Clar83} we know that 
   the Clarke Jacobian of the mapping $H(\cdot,v)$ defined in Lemma \ref{pro-thetak}(d) 
   takes the following form
   \[ 
     \partial_z H(z,v) = \partial\phi(z_1)\times\partial\phi(z_2)\times\cdots\times\partial\phi(z_n)
   \]
  with $\partial\phi(z_i)=\{1\}$ if $v_i=0$ and otherwise 
  $\partial\phi(z_i)=\left\{\begin{array}{cl}
                          \{1\} & {\rm if}\ |z_i|>\lambda v_i,\\
                         \ [0,1] & {\rm if}\ |z_i|=\lambda v_i,\\
                           \{0\} & {\rm if}\ |z_i|<\lambda v_i.
                        \end{array}\right. $
  \end{remark}

  By Lemma \ref{pro-thetak}(a) and (b), we can apply the limited-memory
  BFGS algorithm \cite{NW06} for solving (\ref{dual-PPA1}), but the direction
  yielded by this method may not approximate the Newton direction well
  if the elements in $\wh{\partial}^2\Phi(y^k)$ are badly scaled since
  $\Phi$ is only once continuously differentiable.
  So, we need some Newton steps to bring in the second-order information.
  In view of Lemma \ref{pro-thetak}(c) and (d), we apply the semismooth
  Newton method \cite{QS93} for finding a root of the nonsmooth system
  $\nabla\Phi(y) =0$. To make it possible to solve large-scale problems,
  we use the conjugate gradient (CG) method to yield approximate Newton steps. This leads to
the following   semismooth Newton-CG method.

 \begin{algorithm} \label{Newton-CG}({\bf The semismooth Newton-CG method for (\ref{dual-PPA})})
 \begin{description}
 \item[(S0)] Given $\overline{\eps}>0,\,j_{\max}>0,\,\tau_1,\tau_2\in(0,1),\,
              \varrho\in(0,1)$ and $\mu\in(0,\frac{1}{2})$.
              Choose a starting \hspace*{0.15cm} point $y^0\in\Re^m$ and set $j:=0$.

 \item[(S1)] If $\|\nabla\Phi(y^j)\|\le \overline{\eps}$ or $j>j_{\max}$, then stop.
              Otherwise, go to the next step.

 \item[(S2)] Apply the CG method to seek an approximate solution
              $d^j$ to the linear system
              \ba\label{subprob}
                (V^j+\veps^jI)d=-\nabla\Phi(y^j),
              \ea
             \hspace*{0.1cm} where $V^j\in \wh{\partial}^2\Phi(y^j)$
             with $\wh{\partial}^2\Phi(\cdot)$ given by (\ref{upper-estimate}),
             and $\veps^j:=\tau_1\min\{\tau_2,\|\nabla\Phi(y^j)\|\}$.

 \item[(S3)] Seek the smallest nonnegative integer $l_j$ such that
              the following inequality holds:
              \[
                \Phi(y^{j}+\varrho^{l_j}d^j)\le \Phi(y^j)+\mu\varrho^{l_j}\lla\nabla\Phi(y^j),d^j\rla.
              \]

 \item[(S4)] Set $y^{j+1}:=y^j+\varrho^{l_j}d^j$ and $j:=j+1$, and then go to Step (S.1).
 \end{description}
 \end{algorithm}

  From the definition of $\wh{\partial}^2\Phi(\cdot)$ in Lemma \ref{pro-thetak}(d),
  $V^j$ in Step (S2) of Algorithm \ref{Newton-CG} is positive definite,
  and consequently the search direction $d^j$ is always a descent direction.
  For the global convergence and the rate of local convergence of Algorithm \ref{Newton-CG},
  the interested readers may refer to \cite{ZST10}. Once we have an approximate
  optimal $y^{k}$ of (\ref{dual-PPA}), the approximate optimal solution
  $(u^{k+1},x^{k+1})$ of (\ref{subpro-PPA}) is obtained from the formulas
  \[
    u^{k+1}:=-\beta^{-1}y^{k}\ \ {\rm and}\ \
    x^{k+1}:=S_{\lamd_k}(x^k-\lamd_k A^Ty^k,v).
  \]

  \medskip

  To close this section, we take a look at the selection of $V^j$ in Step (S2)
  for numerical experiments of the next section. By the definition of
  $\wh{\partial}^2\Phi(\cdot)$, $V^j$ takes the form of
  \be\label{Vj}
    V^j = \beta^{-1} I + \lamd_k A D^kA^T,
  \ee
  where $D^k\in \partial_x H(z,v)$ with $z=x^k-\lamd_k A^Ty$. By Remark \ref{Hsubdiff}, 
  $D^k$ is a diagonal matrix, and we select the $i$th diagonal element $D_i^k=1$ if $|z_i|\ge \lambda v_i$ 
  and otherwise $D_i^k=0$.  

 \section{Numerical experiments}\label{Numerical-experiments}

  In this section, we test the performance of Algorithm \ref{Alg1} with the subproblem
 (\ref{subprob11}) solved by the partial PPA. Notice that using the L-BFGS or the semismooth
  Newton-CG alone to solve the subproblem (\ref{subpro-PPA}) of the partial PPA can not yield
  the desired result, since using the L-BFGS alone will not yield good feasibility for
  those difficult problems due to the lack of the second-order information of objective function,
  while using the semismooth Newton-CG alone will meet difficulty for the weighted $l_1$-norm subproblems
  involved in the beginning of Algorithm \ref{Alg1}. In view of this, we develop an exact penalty decomposition
  algorithm with the subproblem (\ref{subprob11}) solved by the partial PPA, for which the subproblems
  (\ref{subpro-PPA}) are solved by combining the L-BFGS with the semismooth Newton-CG.
  The detailed iteration steps of the whole algorithm are described as follows,
 where for any given $\beta_{k},\lamd_{k}\!>0$ and $(x^k,v^k)\!\in\!\Re^n\!\times\!\Re^n$,
 the function $\Phi_k:\Re^m\to\Re$ is defined as
 \[
   \Phi_{k}(y):= b^Ty +\frac{1}{2\beta_{k}}\|y\|^2
   +\frac{1}{2\lamd_{k}}\|S_{\lamd_{k}}(x^{k-1}-\lamd_{k} A^Ty,v^{k-1})\|^2\quad\ \forall y\in\Re^m.
 \]

 \begin{algorithm} \label{Hybrid-Alg1}({\bf Practical exact penalty decomposition method for (\ref{l0})})
 \begin{description}
 \item[(S.0)] Given $\eps,\eps_1>0,\,\omega_1,\omega_2>0$, $\gamma\in(0,1)$, $\sigma\ge 1$ and $\ulamd>0$.
             Choose a sufficiently large $\beta_0$ and suitable $\lambda_0>0$ and $\rho_0>0$.
              Set $(x^0,v^0,y^0)=(0,e,e)$ and $k=0$.

  \item[(S.1)] {\bf While} $\frac{\|Ax^{k}-b\|}{\max\{1,\,\|b\|\}}>\eps_1$
                and $\lamd_{k}>\ulamd$ {\bf do}
               \begin{itemize}
               \item  Set $\lamd_{k+1}=\gamma^{k}\lamd_0$ and $\beta_{k+1}=\beta_{k}$.

                \item  With $y^{k}$ as the starting point, find
                        \(
                          y^{k+1}\!\approx\mathop{\arg\min}_{y\in\Re^{m}} \Phi_{k+1}(y)
                         \)
                        such that $\|\nabla\Phi_{k+1}(y^{k+1})\|\le\omega_1$ by using
                        the L-BFGS algorithm.

              \item  Set $x^{k+1}:=S_{\lamd_{k+1}}(x^{k}-\lamd_{k+1} A^Ty^{k+1},v^{k})$
                      and $v_i^{k+1}:=\left\{\begin{array}{cl}
                                          0 & {\rm if}\ x_i^{k+1}>\rho_k^{-1},\\
                                          1 & {\rm otherwise}.
                                         \end{array}\right.$

               \item  Set $\rho_{k+1}=\sigma\rho_k$ and $k:=k+1$.

              \end{itemize}
              {\bf End}

 \item[(S.2)] {\bf While} $\frac{\|Ax^{k}-b\|}{\max\{1,\,\|b\|\}}>\eps_1$
                or $\lla v^{k},|x^{k}|\rla>\eps$ {\bf do}
             \begin{itemize}
               \item  Set $\lamd_{k+1}=\lamd_{k}$ and $\beta_{k+1}=\beta_{k}$.

               \item  With $y^{k}$ as the starting point, find $y^{k+1}\approx\mathop{\arg\min}_{y\in\Re^{m}}\Phi_{k+1}(y)$
                      such that $\|\nabla\Phi_{k+1}(y^{k+1})\|\le\omega_2$ by using Algorithm \ref{Newton-CG}.

              \item   Set $x^{k+1}:=S_{\lamd_{k+1}}\big(x^{k}-\lamd_{k+1} A^Ty^{k+1},v^{k}\big)$
                       and $v_i^{k+1}:=\left\{\begin{array}{cl}
                                          0 & {\rm if}\ x_i^{k+1}>\rho_k^{-1},\\
                                          1 & {\rm otherwise}.
                                         \end{array}\right.$

             \item Set $\rho_{k+1}= \sigma\rho_k$ and $k:=k+1$.
             \end{itemize}
              {\bf End}
 \end{description}
 \end{algorithm}

 \medskip

  By the choice of starting point $(x^0,v^0,y^0)$, the first step
  of Algorithm \ref{Hybrid-Alg1} is solving
  \[
    \min_{u\in\Re^m,x\in\Re^n}\left\{\|x\|_1+\frac{\beta_0}{2}\|u\|^2+\frac{1}{2\lamd_0}\|x\|^2\!:\ Ax+u=b\right\},
  \]
  whose solution is the minimum-norm solution of the $l_1$-norm minimization problem
  \begin{align}\label{l1-norm}
   \min_{x\in\Re^n} \Big\{\|x\|_1:\ Ax = b\Big\}
  \end{align}
  if $\beta_0$ and $\lamd_0$ are chosen to be sufficiently large (see \cite{MM79}).
  Taking into account that the $l_1$-norm minimization problem is a good convex surrogate for the zero-norm
  minimization problem (\ref{l0}), we should solve the problem $\min_{y\in\Re^{m}}\Phi_{1}(y)$ as
  well as we can. If the initial step can not yield an iterate with good feasibility,
  then we solve the regularized problems
  \be\label{subpro41}
    \min_{x\in\Re^n,u\in\Re^m}\left\{\lla v^k,|x|\rla +\frac{\beta_{k+1}}{2}\|u\|^2+\frac{1}{2\lamd_{k+1}}\|x-x^k\|^2: Ax+u=b\right\}.
  \ee
  with a decreasing sequence $\{\lamd_k\}$ and a nondecreasing sequence $\{\beta_k\}$
  via the L-BFGS. Once a good feasible point is found in Step (S.1),
  Algorithm \ref{Hybrid-Alg1} turns to the second stage, i.e., to solve (\ref{subpro41})
  with nondecreasing sequences $\{\beta_k\}$ and $\{\lamd_k\}$ via Algorithm \ref{Newton-CG}.

  \medskip

  Unless otherwise stated, the parameters involved in Algorithm \ref{Hybrid-Alg1}
  were chosen as:
  \ba\label{para}
    \eps =\frac{10^{-2}}{\max(1,\|b\|)},\,\eps_1=10^{-6},\,\omega_1=10^{-5},\,\omega_2 =10^{-6},\,
    \ulamd=10^{-2},\,\sigma=2,\qquad\nonumber\\
    \beta_0=\max(5\|b\|\times 10^{6},10^{10}),\,\rho_0=\min(1,10/\|b\|),\,
     \lambda_0=\wh{\gamma}\|b\|,\qquad\qquad
  \ea
  where we set $\gamma=0.6$ and $\wh{\gamma}=5$ if $A$ is stored implicitly
  (i.e., $A$ is given in operator form); and otherwise we chose
  $\gamma$ and $\wh{\gamma}$ by the scale of the problem, i.e.,
   \[
    \gamma=\left\{\begin{array}{cl}
                  0.5 & {\rm if}\ \|b\|>10^5\ {\rm or}\ \|b\|\le 5,\\
                  0.8 & {\rm otherwise},
            \end{array}\right. \ {\rm and}\ \
    \wh{\gamma}=\left\{\begin{array}{cl}
                  10 & {\rm if}\ \|b\|>10^5\ {\rm or}\ \|b\|\le 5,\\
                 1.5 & {\rm otherwise}.
            \end{array}\right.
  \]
  We employed the L-BFGS with 5 limited-memory vector-updates and the nonmonotone Armijo line search rule \cite{GLL86}
  to yield an approximate solution to the minimization problem in Step (S.1) of Algorithm \ref{Hybrid-Alg1}. Among others,
  the number of maximum iterations of the L-BFGS was chosen as {\bf 300} for the minimization of $\Phi_1(y)$,
  and {\bf 50} for the minimization of $\Phi_{k}(y)$ with $k\ge 2$.
  The parameters involved in Algorithm \ref{Newton-CG} are set as:
  \ba\label{par-NT}
    \overline{\eps}=10^{-6},\, j_{\rm max}=50,\,
    \tau_1=0.1,\,\tau_2=10^{-4},\,\varrho=0.5,\,\mu=10^{-4}.
  \ea
  In addition, during the testing, if the decrease of gradient is slow in Step (S.1),
  we terminate the L-BFGS in advance and then turn to the solution of the next subproblem.
  Unless otherwise stated, the parameters in {\bf QPDM} and {\bf ISDM} are all set to default values,
  the ``Hybridls" type line search and ``lbfgs" type subspace optimization method are chosen for {\bf FPC\_AS},
   and $\mu=10^{-10},\epsilon=10^{-12}$ and $\epsilon_x=10^{-16}$ are used for {\bf FPC\_AS}.
  All tests described in this section were run in MATLAB R2012(a) under a Windows
  operating system on an Intel Core(TM) i3-2120 3.30GHz CPU with 3GB memory.

  \medskip

  To verify the effectiveness of Algorithm \ref{Hybrid-Alg1}, we compared it with {\bf QPDM} \cite{ZL10},
  {\bf ISDM} \cite{WY10} and ${\bf FPC\_AS}$ on four different sets of problems.
  Since the four solvers return solutions with tiny but nonzero entries that can be regarded as zero,
  we use {\bf nnzx} to denote the number of nonzeros in $x$ which we estimate
  as in \cite{BF08} by the minimum cardinality of a subset of the components
  of $x$ that account for $99.9\%$ of $\|x\|_1$; i.e.,
  \[
    {\bf nnzx}:=\min\left\{\kappa:\ts\sum_{i=1}^\kappa|x|_i^{\downarrow} \ge 0.999\|x\|_1\right\}.
  \]
  Suppose that the exact sparsest solution $x^*$ is known. We also compare
  the support of $x^f$ with that of $x^*$, where $x^{f}$ is the final iterate yielded
  by the above four solvers. To this end, we first remove tiny entries of $x^{f}$
  by setting all of its entries with a magnitude smaller than $0.1|x^*|_{\rm snz}$
  to zero, where $|x^*|_{\rm snz}$ is the smallest nonzero component of $|x^*|$,
  and then compute the quantities ``{\bf sgn}", ``{\bf miss}" and ``{\bf over}", where
  \[
    {\bf sgn}:=\big|\{i\ |\ x_i^{f}x_i^*<0\}\big|,\
    {\bf miss}:=\big|\{i\ |\ x_i^{f}=0,x_i^*\ne 0\}\big|,\,
    {\bf over}:=\big|\{i\ |\ x_i^{f}\ne 0,x_i^*=0\}\big|.
  \]


  \subsection{Recoverability for some ``pathological" problems} \label{Subsec5.1}

   We tested Algorithm \ref{Hybrid-Alg1}, {\bf FPC\_AS}, {\bf ISDM} and {\bf QPDM} on a set of small-scale,
  pathological problems described in Table \ref{tab1}. The first test set includes four problems
  Caltech Test 1, \ldots, Caltech Test 4 given by Cand\`{e}s and Becker, which, as mentioned in \cite{WYGZ10},
  are pathological because the magnitudes of the nonzero entries of the exact solution $x^*$ lies in a large range.
  Such pathological problems are exaggerations of a large number of realistic problems
  in which the signals have both large and small entries. The second test set includes
  six problems Ameth6Xmeth20-Ameth6Xmeth24 and Ameth6Xmeth6 from \cite{WYGZ10},
  which are difficult since the number of nonzero entries in their solutions is close to
  the limit where the zero-norm problem (\ref{l0}) is equivalent to the $l_1$-norm problem.
 \begin{table}[htbp]
 \caption{\label{tab1} Description of some pathological problems}
 \begin{tabular}{|c|c|c|c|c|c|}
  \hline
   \ {\bf ID} & Name & $n$ & $m$ & $K$ & {\rm (Magnitude, num. of entries on this level)} \\
  \hline
    \  1 & {\rm CaltechTest1} & 512 & 128 & 38 & $(10^5,33),(1,5)$\\
  \hline
   \  2 & {\rm CaltechTest2} & 512 & 128 & 37 & $(10^5,32),(1,5)$\\
   \hline
   \  3 & {\rm CaltechTest3} & 512 & 128 & 32 & $(10^5,31),(10^{-6},1)$\\
  \hline
   \  4 & {\rm CaltechTest4} & 512 & 102 & 26 & $(10^4,13),(1,12),(10^{-2},1)$\\
   \hline
   \  5 & {\rm Ameth6Xmeth20} & 1024 & 512 & 150 & $(1,150)$\\
   \hline
   \  6 & {\rm Ameth6Xmeth21} & 1024 & 512 & 151 & $(1,150)$\\
   \hline
   \  7 & {\rm Ameth6Xmeth22} & 1024 & 512 & 152 & $(1,150)$\\
   \hline
   \  8 & {\rm Ameth6Xmeth23} & 1024 & 512 & 153 & $(1,150)$\\
   \hline
   \  9 & {\rm Ameth6Xmeth24} & 1024 & 512 & 154 & $(1,150)$\\
   \hline
   \  10 & {\rm Ameth6Xmeth6} & 1024 & 512 & 154 & $(1,150)$\\
  \hline
  \end{tabular}
  \end{table}

 \begin{table}[htbp]
 \begin{center}
 \caption{\label{tab2} Numerical results of four solvers for the pathological problems}
 \begin{tabular}{|c|c|c|c|c|c|c|c|}
  \hline
   \ {\bf ID} & {\bf Solver} & {\bf time(s)} & {\bf Relerr} & {\bf Res} & {\bf nMat} & {\bf nnzx}
     & ({\bf sgn}, {\bf miss}, {\bf over})\\
  \hline
    \    & {\bf Algorithm} \ref{Hybrid-Alg1} & 0.39 & 5.16e-12 & 8.87e-9 & 1057 & 33 & $(0,0,0)$\\
    \  1 & {\bf FPC\_AS} & 0.47 & 4.98e-12 & 4.37e-8 & 437 & 33 & $(0,0,0)$\\
    \    & {\bf ISDM} & 0.47 & 4.52e-6 & 9.97e-1 &--  & 33 & $(0,5,0)$\\
     \   & {\bf QPDM} & 0.15 & 2.07e-0 & 1.53e-9 &--  & 125 & $(0,28,118)$\\
    \hline
    \    & {\bf Algorithm} \ref{Hybrid-Alg1} & 0.11 & 8.25e-14 & 8.44e-9 & 1060 & 32 & $(0,0,0)$\\
    \  2 & {\bf FPC\_AS} & 0.12 & 1.86e-13 & 5.58e-8 & 357 & 32 & $(0,0,0)$\\
    \    & {\bf ISDM}  & 0.37 & 4.27e-6 & 9.34e-1 & -- & 32 & $(0,5,0)$\\
     \   & {\bf QPDM} & 0.01 & 2.18e-0 & 1.17e-9 & -- & 123 & $(0,28,119)$\\
    \hline
    \    & {\bf Algorithm} \ref{Hybrid-Alg1} & 0.16 & 4.56e-9 & 2.09e-14 & 1199 & 31 & $(0,0,0)$\\
    \  3 & {\bf FPC\_AS} & 0.06 & 1.15e-9 & 1.61e-9 & 247 & 31 & $(0,0,0)$\\
    \    & {\bf ISDM}  & 0.42  & 9.78e-7 & 4.67e-7 & -- & 31 & $(0,1,0)$\\
         & {\bf QPDM} & 0.03 & 9.78e-7 & 4.67e-7 &--  & 31 & $(0,1,0)$\\
    \hline
    \    & {\bf Algorithm} \ref{Hybrid-Alg1} & 0.16 & 3.10e-7 & 4.05e-3 & 985 & 13 & $(0,1,0)$\\
    \ 4  & {\bf FPC\_AS} & 0.17 & 4.52e-13 & 7.51e-9 & 572 & 13 & $(0,0,0)$\\
    \    & {\bf ISDM}  & 0.34  & 9.96e-5 & 1.38e-0 &--  & 13 & $(0,12,1)$\\
         & {\bf QPDM} & 0.02 & 2.10e-0 & 6.00e-11 & -- & 13 & $(0,21,97)$\\
    \hline
    \    & {\bf Algorithm} \ref{Hybrid-Alg1} & 0.47 & 4.93e-14 & 5.69e-13 & 1098 & 150 & $(0,0,0)$\\
    \ 5  & {\bf FPC\_AS} & 0.25 & 6.80e-10 & 4.01e-9 & 412 & 150 & $(0,0,0)$\\
    \    & {\bf ISDM}  & 3.17  & 6.67e-1 & 3.72e-1 &--  & 464 & $(0,15,185)$\\
         & {\bf QPDM} & 1.58  & 8.65e-1 & 4.38e-1 & -- & 492 & $(0,28,287)$\\
     \hline
    \    & {\bf Algorithm} \ref{Hybrid-Alg1} & 0.36 & 4.91e-14 & 5.68e-13 & 730 & 151 & $(0,0,0)$\\
    \ 6  & {\bf FPC\_AS} & 0.29 & 6.96e-10 & 4.11e-9 & 408 & 151 & $(0,0,0)$\\
    \    & {\bf ISDM}  & 3.38  & 4.92e-14 & 5.80e-14 & -- & 151 & $(0,0,0)$\\
         & {\bf QPDM} & 0.92  & 6.58e-1 & 4.92e-1 & -- & 480 & $(0,16,211)$\\
    \hline
    \    & {\bf Algorithm} \ref{Hybrid-Alg1} & 0.41 & 4.91e-14 & 5.69e-13 & 910 & 152 & $(0,0,0)$\\
    \ 7  & {\bf FPC\_AS} & 0.36 & 8.10e-10 & 4.81e-9 & 461 & 152 & $(0,0,0)$\\
    \    & {\bf ISDM}  & 3.29  & 5.02e-14 & 5.80e-13 &--  & 152 & $(0,0,0)$\\
         & {\bf QPDM} & 1.42  & 6.89e-1 & 5.00e-1 & -- & 481 & $(0,21,222)$\\
     \hline
     \    & {\bf Algorithm} \ref{Hybrid-Alg1} & 0.34 & 4.95e-14 & 5.68e-13 & 809 & 153 & $(0,0,0)$\\
    \ 8  & {\bf FPC\_AS} & 0.34 & 9.19e-10 & 5.44e-9 & 578 & 153 & $(0,0,0)$\\
    \    & {\bf ISDM}  & 2.96  & 5.51e-14 & 5.88e-13 & -- & 153 & $(0,0,0)$\\
         & {\bf QPDM} & 1.36  & 2.06e-0 & 2.50e-1 & -- & 494 & $(0,37,361)$\\
     \hline
     \    & {\bf Algorithm} \ref{Hybrid-Alg1} & 0.37 & 4.94e-14 & 5.67e-13 & 826 & 154 & $(0,0,0)$\\
    \ 9   & {\bf FPC\_AS} & 0.41 & 9.41e-10 & 5.57e-9 & 572 & 154 & $(0,0,0)$\\
    \    & {\bf ISDM}  & 3.81  & 4.70e-14 & 5.79e-13 & -- & 154 & $(0,0,0)$\\
         & {\bf QPDM} & 1.47  & 2.71e-0 & 2.31e-1 &--  & 496 & $(0,42,374)$\\
    \hline
     \    & {\bf Algorithm} \ref{Hybrid-Alg1} & 0.61 & 4.94e-14 & 5.67e-8 & 1419 & 154 & $(0,0,0)$\\
    \ 10  & {\bf FPC\_AS} & 0.31 & 2.59e-13 & 1.57e-7 & 577 & 154 & $(0,0,0)$\\
    \     & {\bf ISDM} & 3.28 & 4.70e-14 & 5.79e-8 & -- & 154 & $(0,0,0)$\\
     \    & {\bf QPDM} & 0.22 & 3.01e-0 & 2.22e+4 & -- & 499 & $(0,81,420)$\\
    \hline
  \end{tabular}
  \end{center}
  \end{table}

  The numerical results of four solvers are reported in Table \ref{tab2},
  where {\bf nMat} means the total number of matrix-vector products involving $A$ and $A^T$,
  {\bf Time} means the computing time in seconds, {\bf Res} denotes the $l_2$-norm of
  recovered residual, i.e., ${\bf Res}=\|Ax^{f}-b\|$, and {\bf Relerr} means
  the relative error between the recovered solution $x^f$ and the true solution $x^*$, i.e.,
  \(
    {\bf Relerr} ={\|x^{f}-x^*\|}/{\|x^*\|}.
  \)
  Since {\bf ISDM} and {\bf QPDM} do not record the number of matrix-vector products
  involving $A$ and $A^T$, we mark {\bf nMat} as ``--".

  \medskip

  Table \ref{tab2} shows that among the four solvers, {\bf QPDM} has the worst performance
  and can not recovery any one of these problems, {\bf ISDM} requires the most computing time
  and yields solutions with incorrect {\bf miss} for the first test set, and Algorithm \ref{Hybrid-Alg1}
  and {\bf FPC\_AS} have comparable performance in terms of recoverability and computing time.

 \subsection{Sparse signal recovery from noiseless measurements}\label{Subsec5.2}

  In this subsection we compare the performance of Algorithm \ref{Hybrid-Alg1} with
  that of {\bf FPC\_AS}, {\bf ISDM} and {\bf QPDM} for compressed sensing reconstruction
  on randomly generated problems. Given the dimension $n$ of a signal,
  the number of observations $m$ and the number of nonzeros $K$, we generated a random matrix
  $A\in\Re^{m\times n}$ and a random $x^*\in\Re^n$ in the same way as in \cite{WYGZ10}.
  Specifically, we generated a matrix by one of the following types:
  \begin{description}
   \item[{\bf Type 1}:] Gaussian matrix whose elements are generated independently and
                        identically \hspace*{0.6cm} distributed from the normal distribution
                        $N(0,1)$;

   \item[{\bf Type 2}:] Orthogonalized Gaussian matrix whose rows are orthogonalized using
                        a QR \hspace*{0.6cm} decomposition;

   \item[{\bf Type 3}:] Bernoulli matrix whose elements are $\pm 1$ independently with equal
                        probability;

   \item[{\bf Type 4}:] Hadamard matrix $H$, which is a matrix of $\pm 1$ whose columns are orthogonal;

   \item[{\bf Type 5}:] Discrete cosine transform (DCT) matrix;
  \end{description}
   and then randomly selected $m$ rows from this matrix to construct the matrix $A$.
   Similar to \cite{WYGZ10}, we also scaled the matrix $A$ constructed from
   matrices of types $1$, $3$, and $4$ by the largest eigenvalue of $AA^T$.
   In order to generate the signal $x^*$,  we first generated the support
   by randomly selecting $K$ indexed between $1$ and $n$, and then assigned
   a value to $x_i^*$ for each $i$ in the support by one of the following six methods:
   \begin{description}
   \item[{\bf Type 1}:] A normally distributed random variable (Gaussian signal);

   \item[{\bf Type 2}:] A uniformly distributed random variable in $(-1,1)$;

   \item[{\bf Type 3}:] One (zero-one signal);

   \item[{\bf Type 4}:] The sign of a normally distributed random variable;



   \item[{\bf Type 5}:] A signal $x$ with power-law decaying entries (known as compressible
                        sparse \hspace*{0.8cm} signals) whose components satisfy $|x_i|\le c_x i^{-p}$,
                        where $c_x=10^5$ and $p=1.5$;

   \item[{\bf Type 6}:] A signal $x$ with exponential decaying entries whose components satisfy
                         \[
                           |x_i|\le c_x e^{-pi}\ \ {\rm with}\ \ c_x=1\ {\rm and}\ p=0.005.
                         \]
  \end{description}
  Finally, the observation $b$ was computed as $b=Ax^*$. The matrices of types $1,2,3$
  and $4$ were stored explicitly, and the matrices of type $5$ were stored implicitly.
  Unless otherwise stated, in the sequel, we call a signal recovered successfully by
  a solver if the relative error between the solution $x^f$ generated and the original signal
  $x^*$ is less than $5\times{\bf 10^{-7}}$.

\begin{figure}[htbp]
 \begin{center}
 {\includegraphics[width=0.95\textwidth]{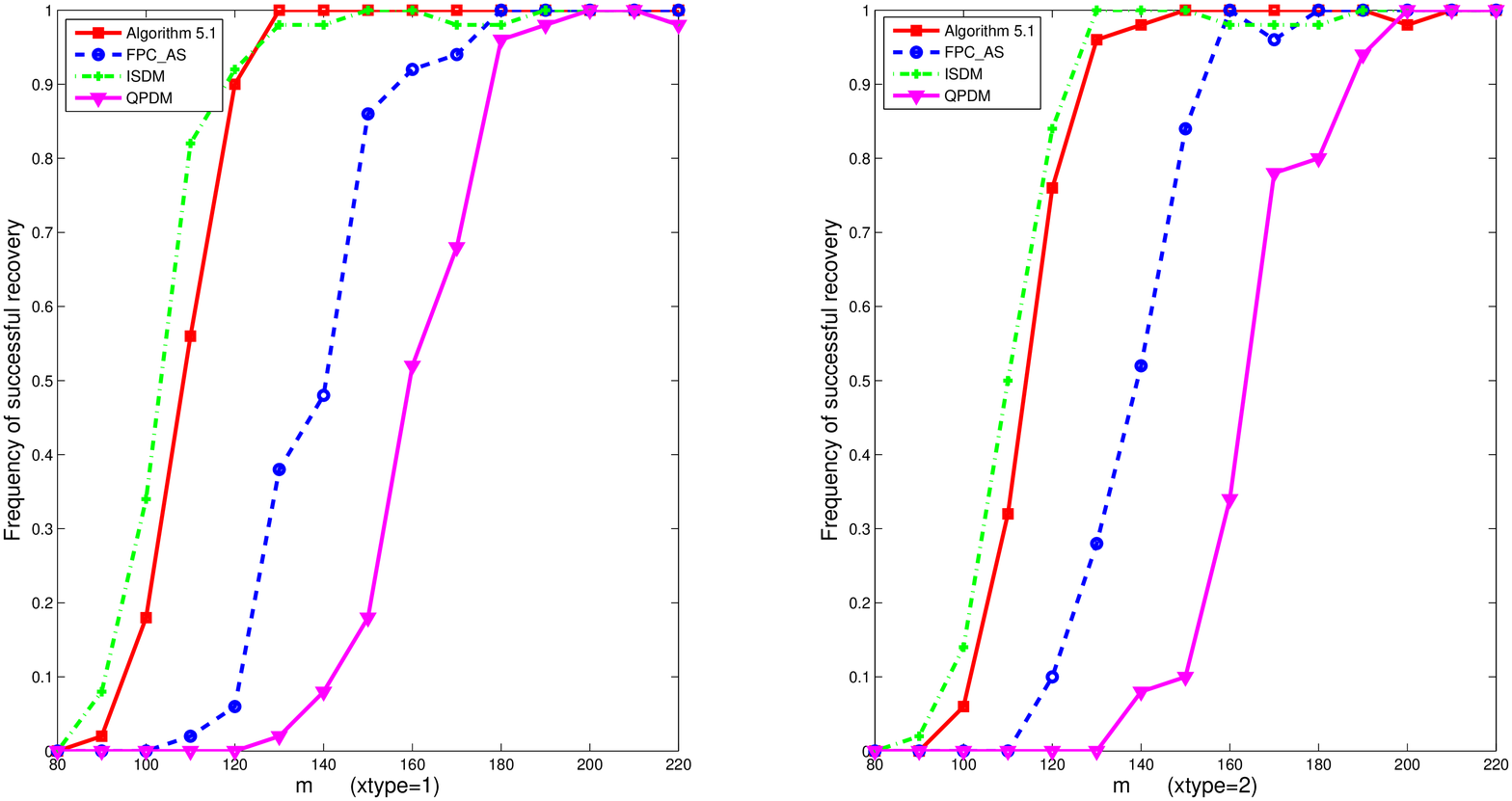}}\\
 {\includegraphics[width=0.95\textwidth]{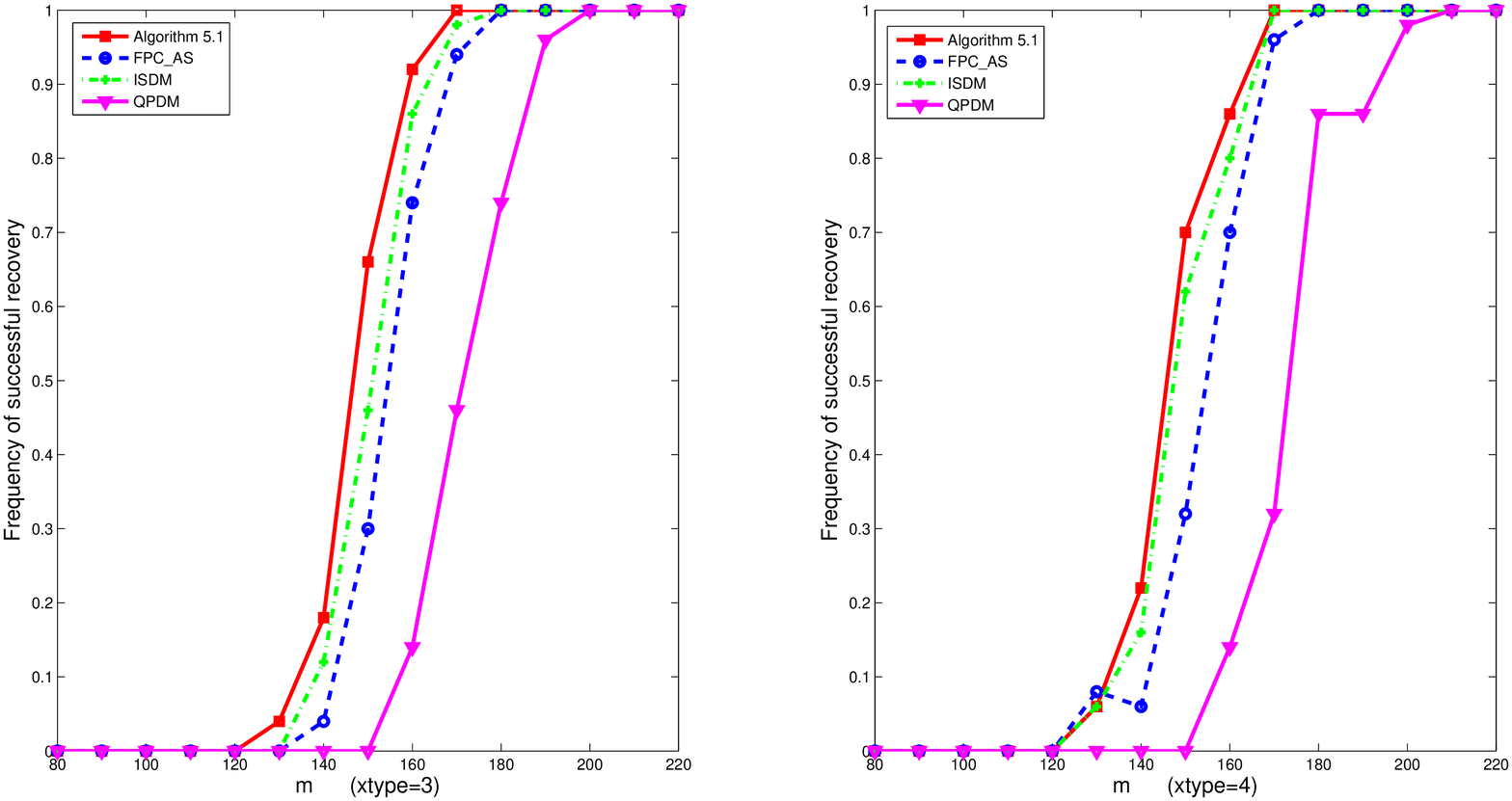}}\\
 {\includegraphics[width=0.95\textwidth]{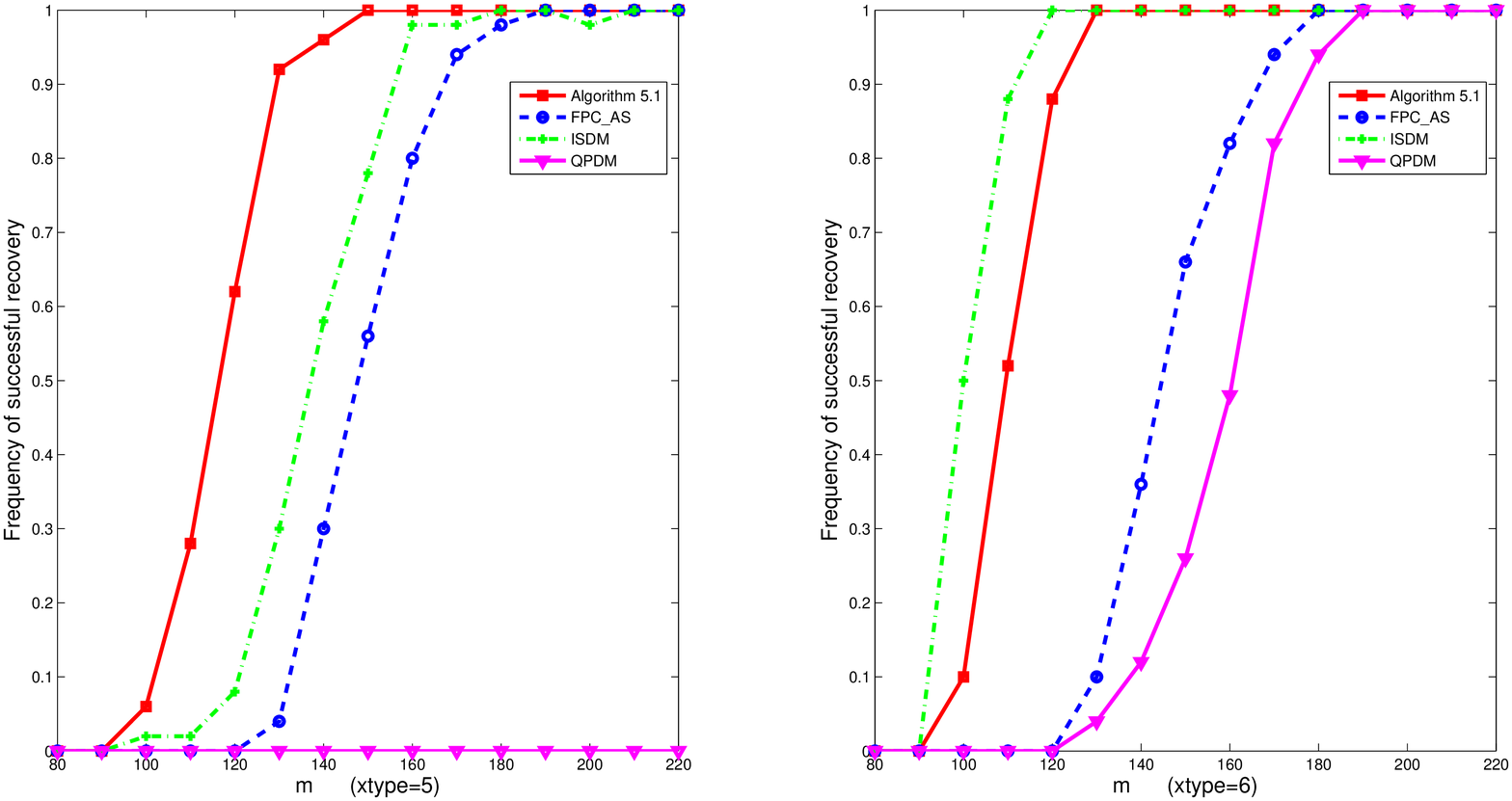}}
 \caption{\small Frequency of successful recovery for four solvers (Atype$=1$)}
 \label{fig1}
 \end{center}
 \end{figure}

\begin{figure}[htbp]
 \begin{center}
 {\includegraphics[width=0.95\textwidth]{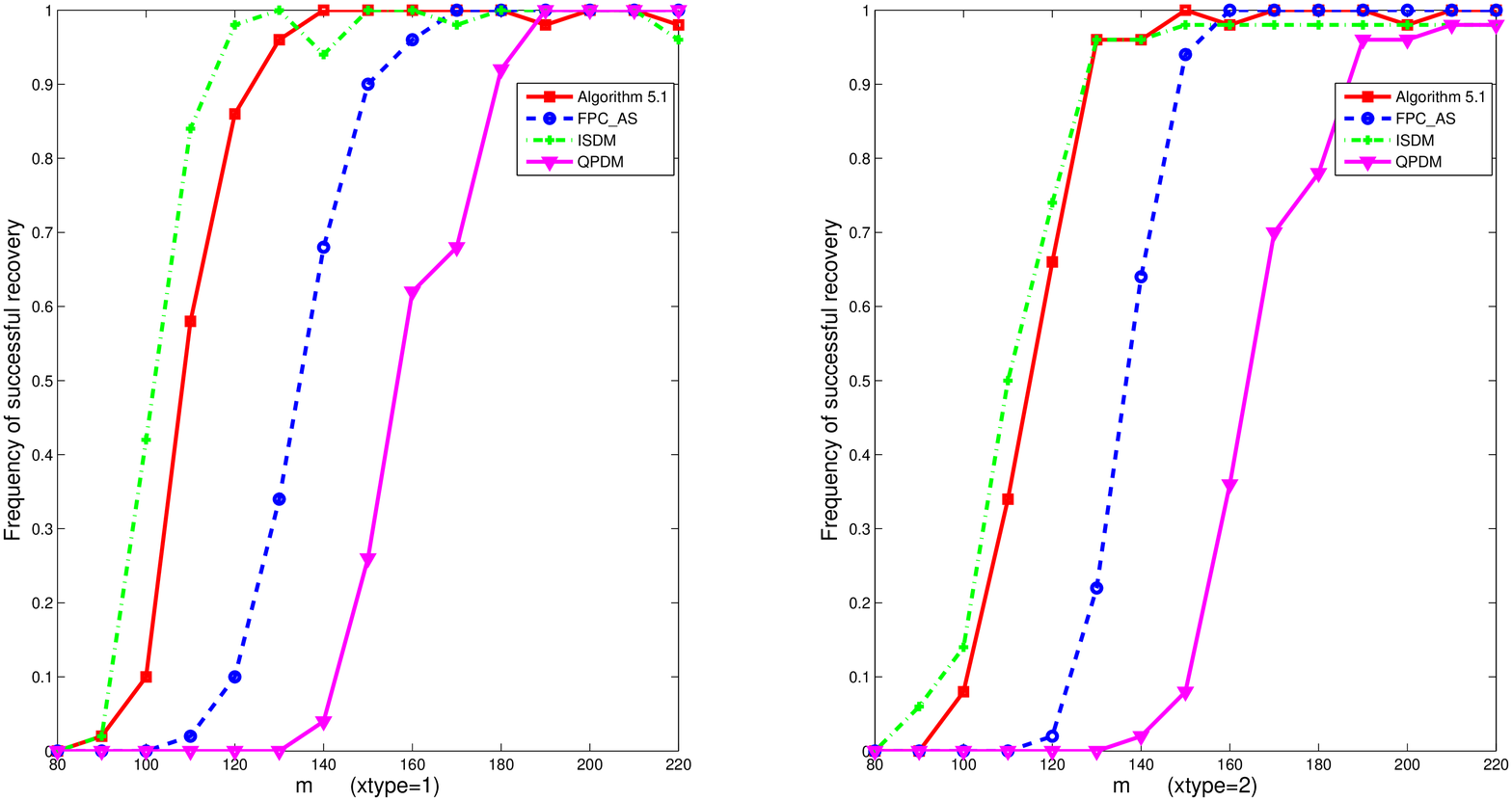}}\\
 {\includegraphics[width=0.95\textwidth]{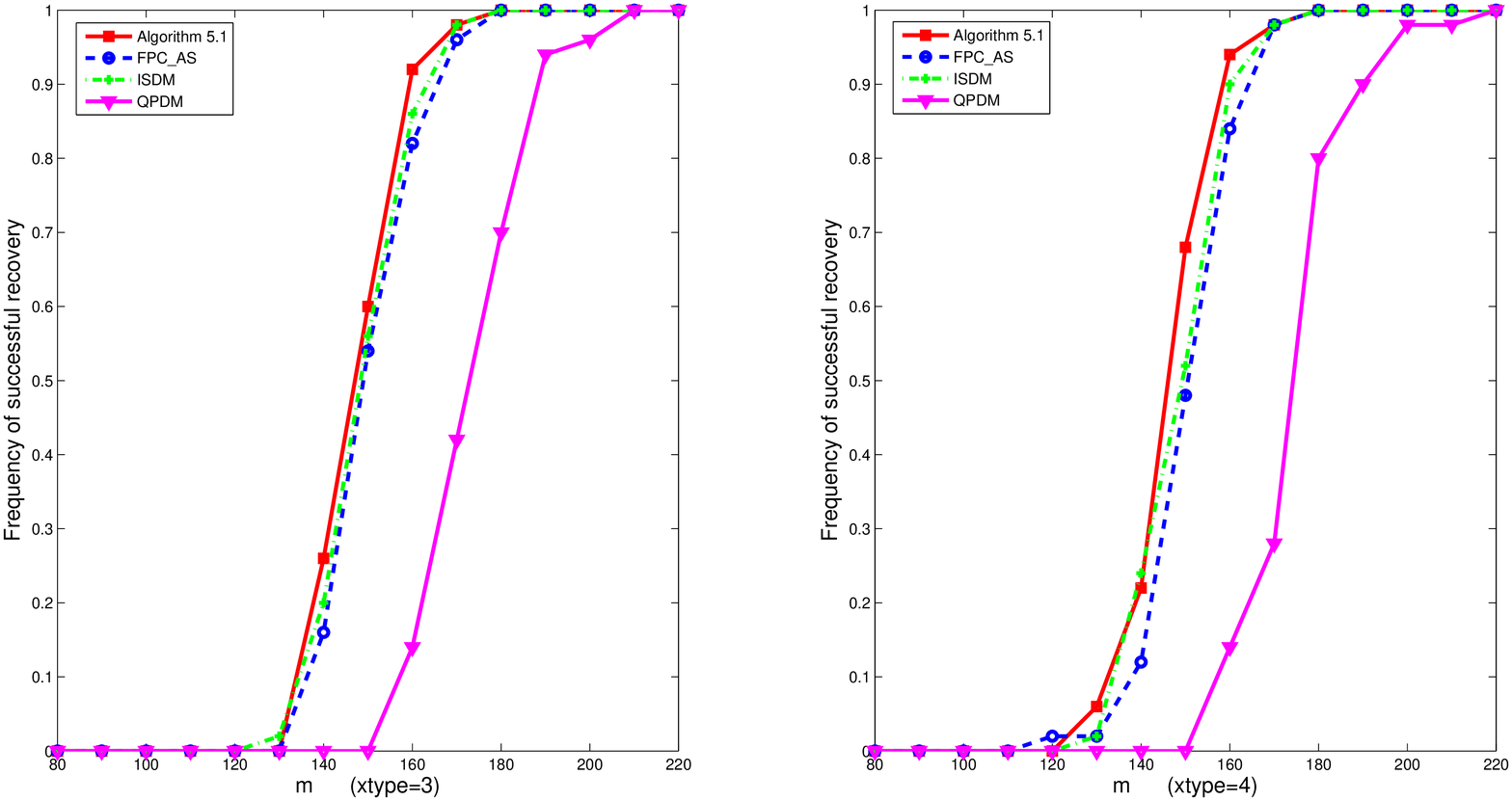}}\\
 {\includegraphics[width=0.95\textwidth]{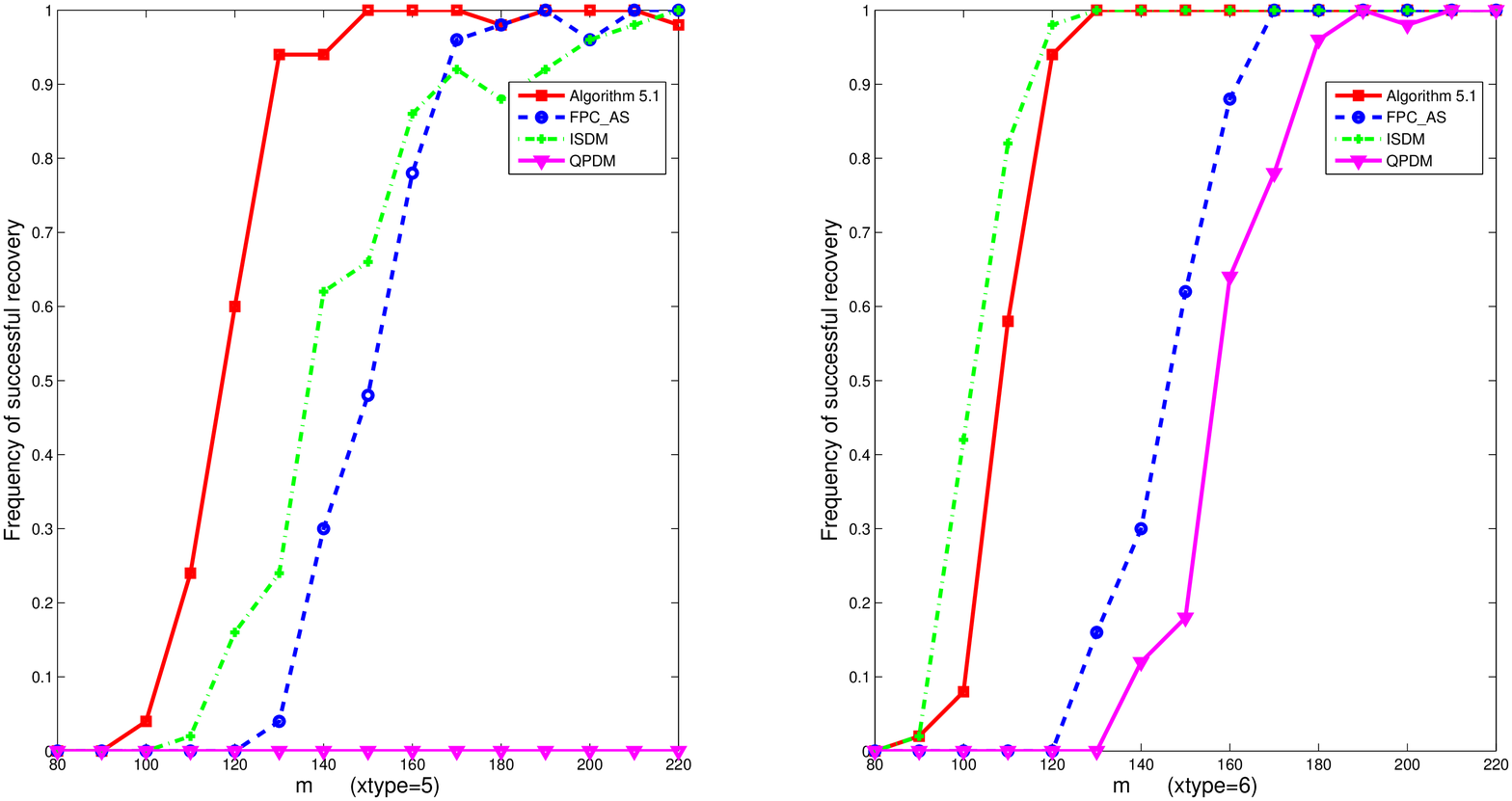}}
 \caption{\small Frequency of successful recovery for four solvers (Atype$=2$)}
 \label{fig2}
 \end{center}
 \end{figure}

 \begin{figure}[htbp]
 \begin{center}
 {\includegraphics[width=0.95\textwidth]{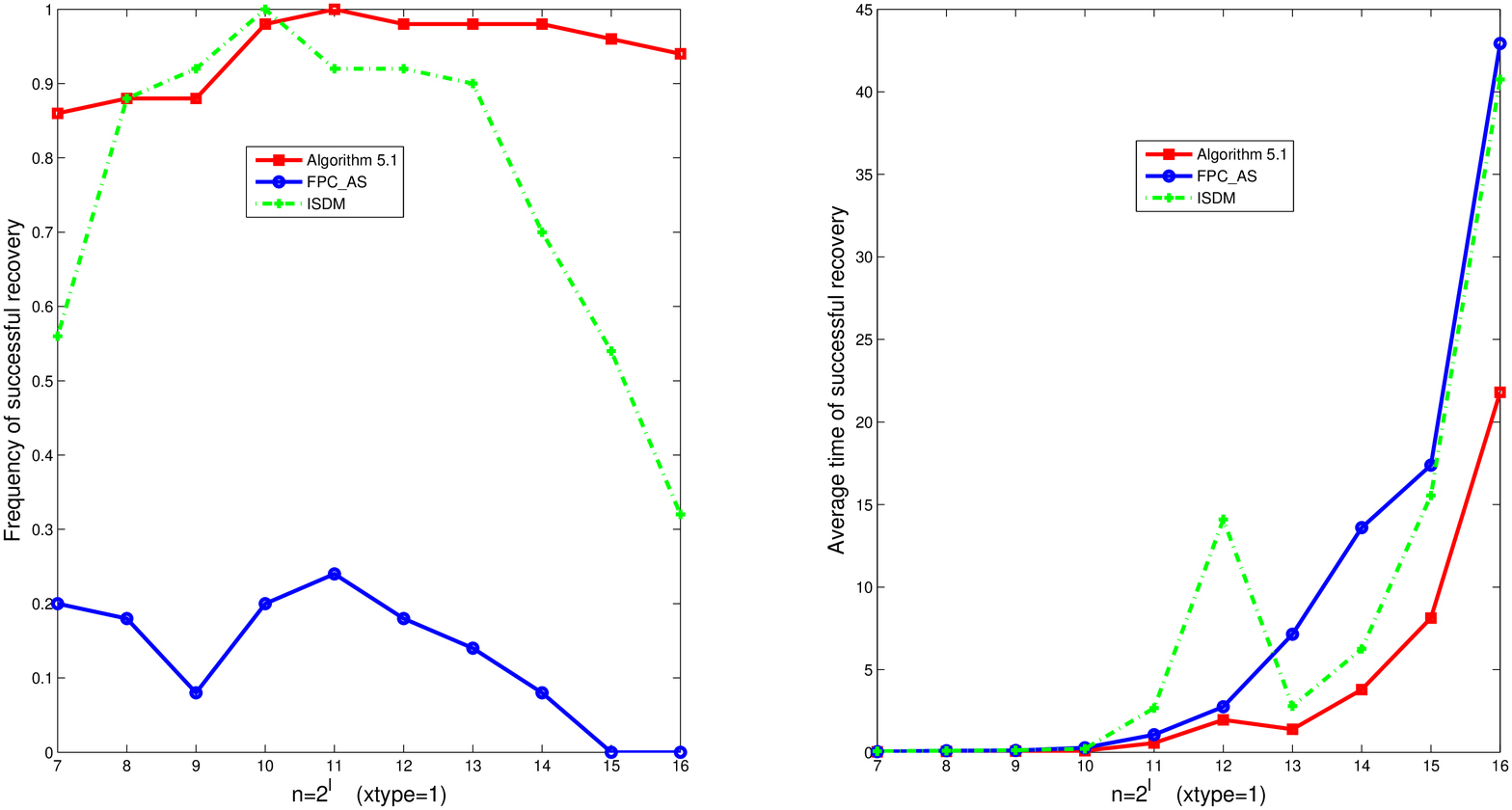}}\\
 {\includegraphics[width=0.95\textwidth]{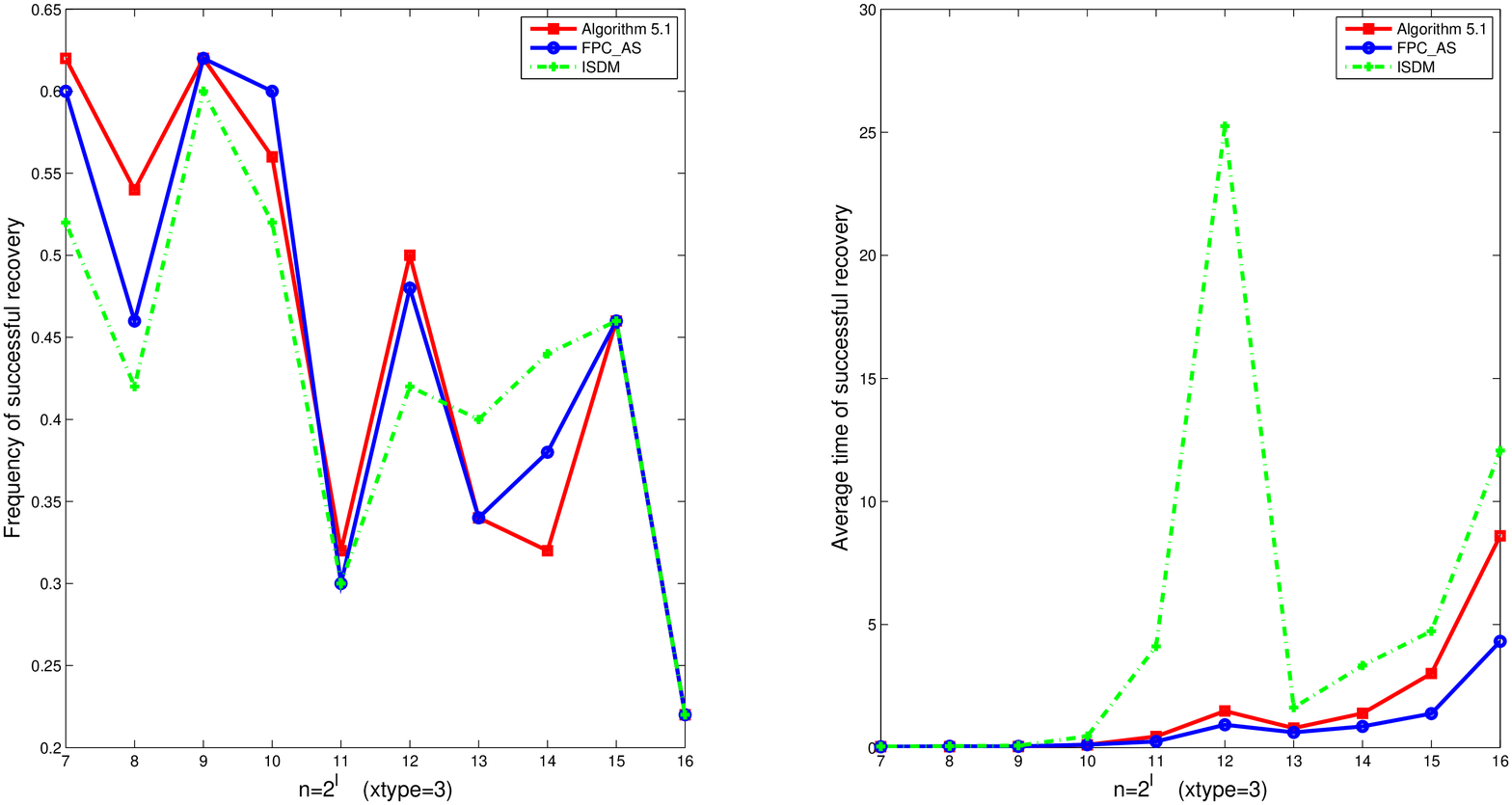}}\\
 {\includegraphics[width=0.95\textwidth]{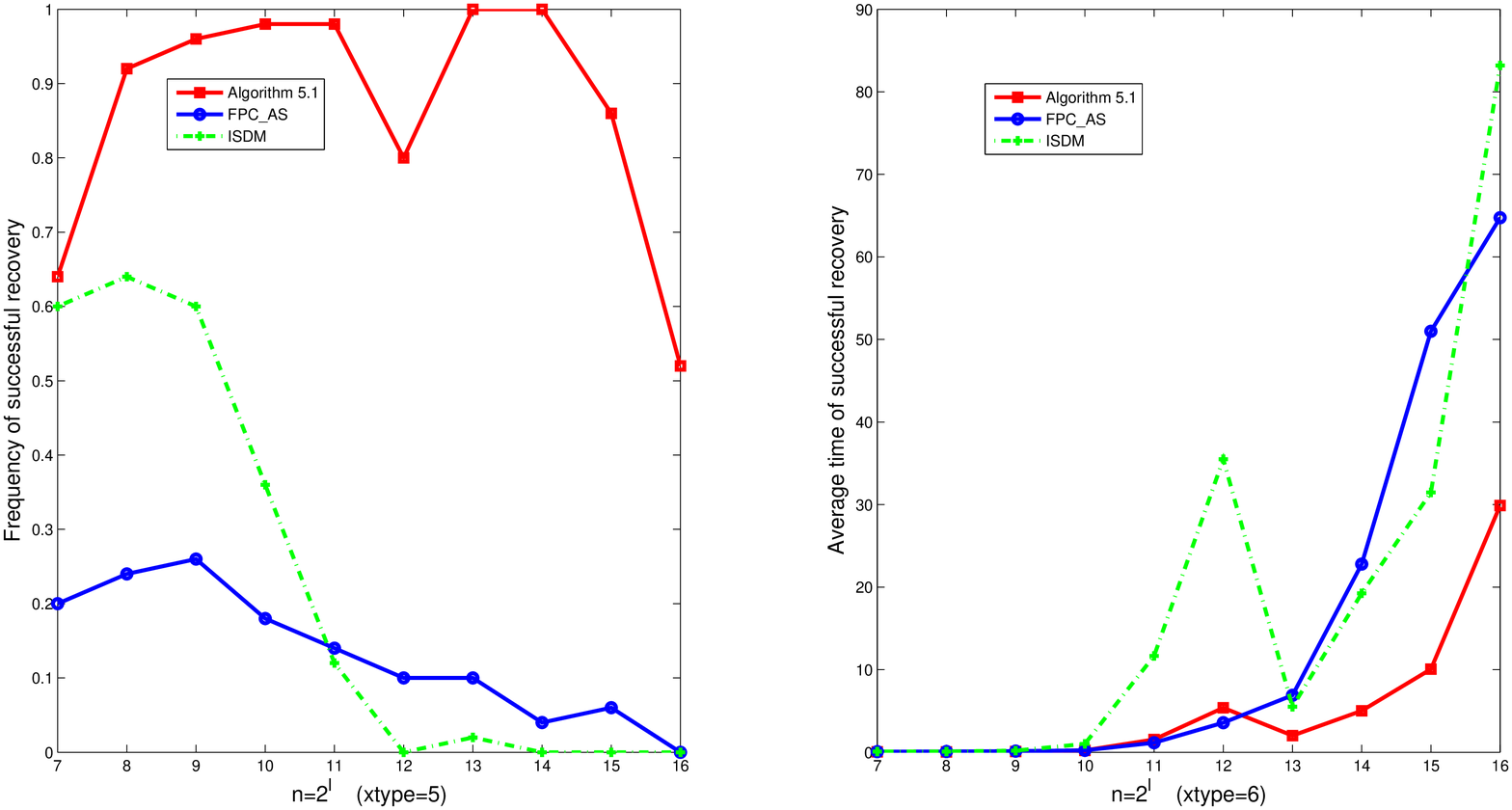}}
 \caption{\small Frequency and time of successful recovery for three solvers (Atype$=5$)}
 \label{fig3}
 \end{center}
 \end{figure}

  \medskip

  We first took the matrix of type $1$ for example to test the influence of
  the number of measurements $m$ on the recoverability of four solvers for
  different types of signals. For each type of signal, we considered
  the dimension $n=600$ and the number of nonzeros $K=40$ and took the number of
  measurements $m\in\{80,90,100,\cdots,220\}$. For each $m$, we generated {\bf 50} problems randomly,
  and tested the frequency of successful recovery for each solver.
  The curves of Figure \ref{fig1} depict how the recoverability of four solvers vary with
  the number of measurements for different types of signals.

  \medskip

  Figure \ref{fig1} shows that among the four solvers, {\bf QPDM} has the worst recoverability
  for all six different types of signals, {\bf ISDM} has a little better recoverability than
  Algorithm \ref{Hybrid-Alg1} for the signals of types $1$, $2$ and $6$, which are much better than
  that of {\bf FPC\_AS}, and Algorithm \ref{Hybrid-Alg1}, {\bf FPC\_AS} and {\bf ISDM} have comparable
  recoverability for the signals of types $3$ and $4$. For the signals of type $5$, Algorithm \ref{Hybrid-Alg1}
  has much better recoverability than {\bf ISDM} and {\bf FPC\_AS}. After further testing, we found that
  for other types of $A$, the four solvers display the similar performance as in Figure \ref{fig1}
  for the six kinds of signals (see Figure \ref{fig2} for type $2$), and {\bf ISDM} requires
  the most computing time among the solvers.

  \medskip

  Then we took the matrices of type $5$ for example to show how the performance of Algorithm \ref{Hybrid-Alg1},
  {\bf FPC\_AS} and {\bf ISDM} scales with the size of the problem. Since Figure \ref{fig1}-\ref{fig2} illustrates
  that the three solvers have the similar performance for the signals of types $1$, $2$ and $6$,
  and the similar performance for the signals of type $3$ and $4$, we compared their performance only
  for the signals of types $1$, $3$ and $5$. For each type of $x^*$, we generated {\bf 50} problems randomly
  for each $n\in\{2^7,2^8,\ldots,2^{16}\}$ and tested the frequency and the average time of successful recovery
  for the three solvers, where $m={\bf round}(n/6)$ for
  the signals of type $1$, $m={\bf round}(n/3)$ for the signals of type $3$, and $m={\bf round}(n/4)$ for
  the signals of type $5$, and the number of nonzeros $K$ was set to ${\bf round}(0.3m)$.
  The curves of Figure \ref{fig3} depict how the recoverability and the average time of successful recovery
  vary with the size of the problem. When there is no signal recovered successfully,
  we replace the average time of successful recovery by the average computing time of $50$ problems.

  \medskip

  Figure \ref{fig3} shows that for the signals of types $1$ and $5$, Algorithm \ref{Hybrid-Alg1} has much
  higher recoverability than {\bf FPC\_AS} and {\bf ISDM} and requires the less recovery time;
  for the signals of type $3$, the recoverability and the recovery time of three solvers are comparable.

  \medskip

  From Figure \ref{fig1}-\ref{fig3}, we conclude that for all types of matrices considered,
  Algorithm \ref{Hybrid-Alg1} has comparable even better recoverability than {\bf ISDM}
  for the six types of signals above, and requires less computing time than {\bf ISDM};
  Algorithm \ref{Hybrid-Alg1} has better recoverability than {\bf FPC\_AS} and needs comparable even
  less computing time than {\bf FPC\_AS}; and {\bf QPDM} has the worst recoverability for all types of signals.
  In view of this, we did not compare Algorithm \ref{Hybrid-Alg1} with {\bf QPDM}
  in the subsequent numerical experiments.

 \subsection{Sparse signal recovery from noisy measurements}\label{Subsec5.3}

  Since problems in practice are usually corrupted by noise, in this subsection
  we test the recoverability of Algorithm \ref{Hybrid-Alg1} on the same matrices
  and signals as in Subsection \ref{Subsec5.2} but with Gaussian noise,
  and compare its performance with that of {\bf FPC\_AS} and {\bf ISDM}.
  Specifically, we let
  \(
     b = Ax^* + \theta{\xi}/{\|\xi\|},
  \)
  where $\xi$ is a vector whose components are independently and identically distributed as $N(0,1)$,
  and $\theta>0$ is a given constant to denote the noise level. During the testing,
  we always set $\theta$ to $0.01$, and chose the parameters of Algorithm \ref{Hybrid-Alg1}
  as in (\ref{para}) and (\ref{par-NT}) except that $\eps=1,\ \eps_1=\frac{0.01\theta}{\max(1,\|b\|)}$,
  \[
     \gamma=\left\{\begin{array}{cl}
                  0.5 & {\rm if}\ \|b\|\ge 10^2,\\
                  0.8 & {\rm otherwise},
            \end{array}\right.\ \ \wh{\gamma}=\left\{\begin{array}{cl}
                  1 & {\rm if}\ \|b\|\ge 10^2,\\
                 10 & {\rm otherwise}.
            \end{array}\right.\ \ {\rm and}\ \ j_{\rm max}=5.
  \]

 \begin{figure}[htbp]
 \begin{center}
 {\includegraphics[width=0.95\textwidth]{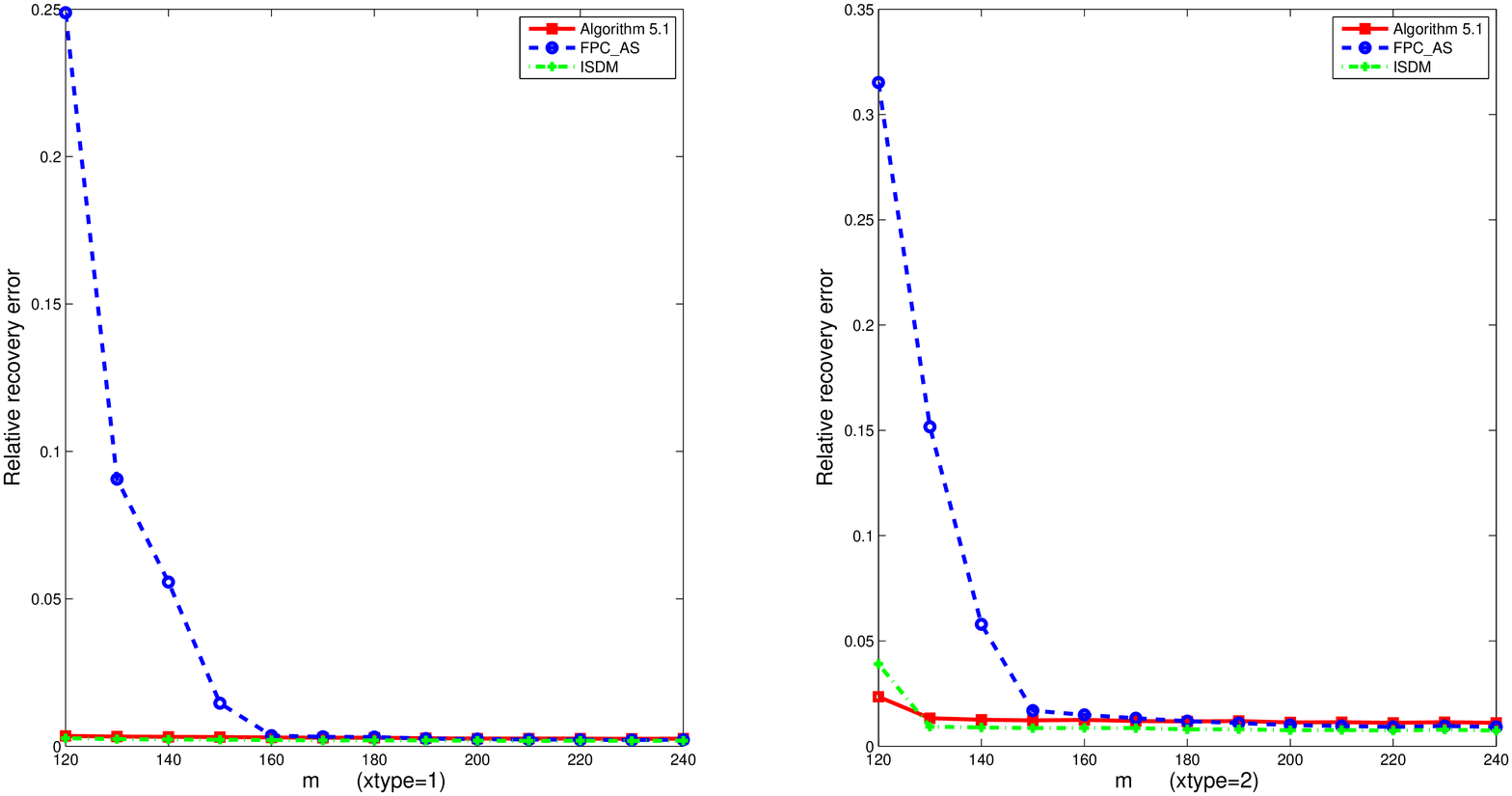}}\\
 {\includegraphics[width=0.95\textwidth]{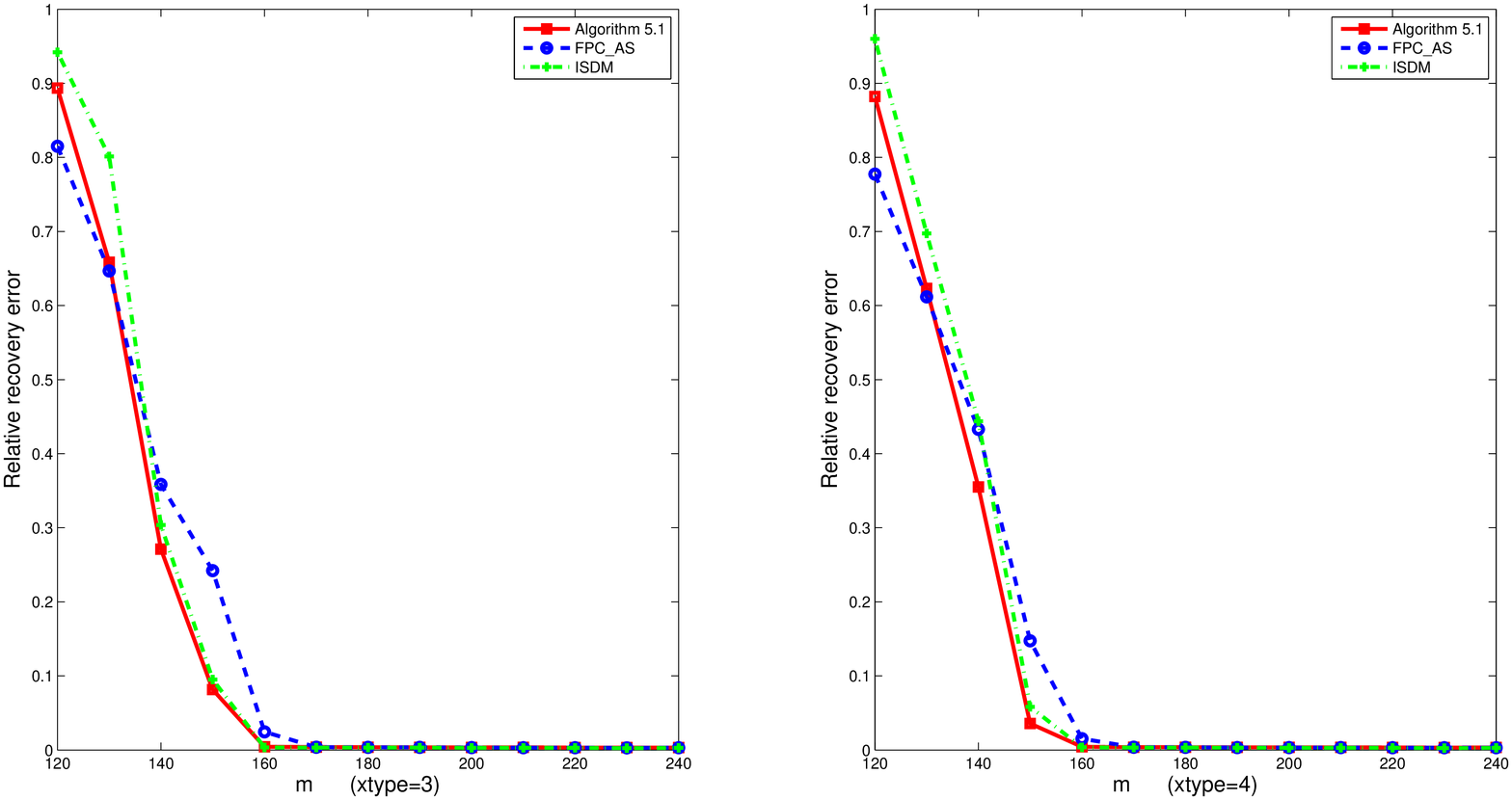}}\\
 {\includegraphics[width=0.95\textwidth]{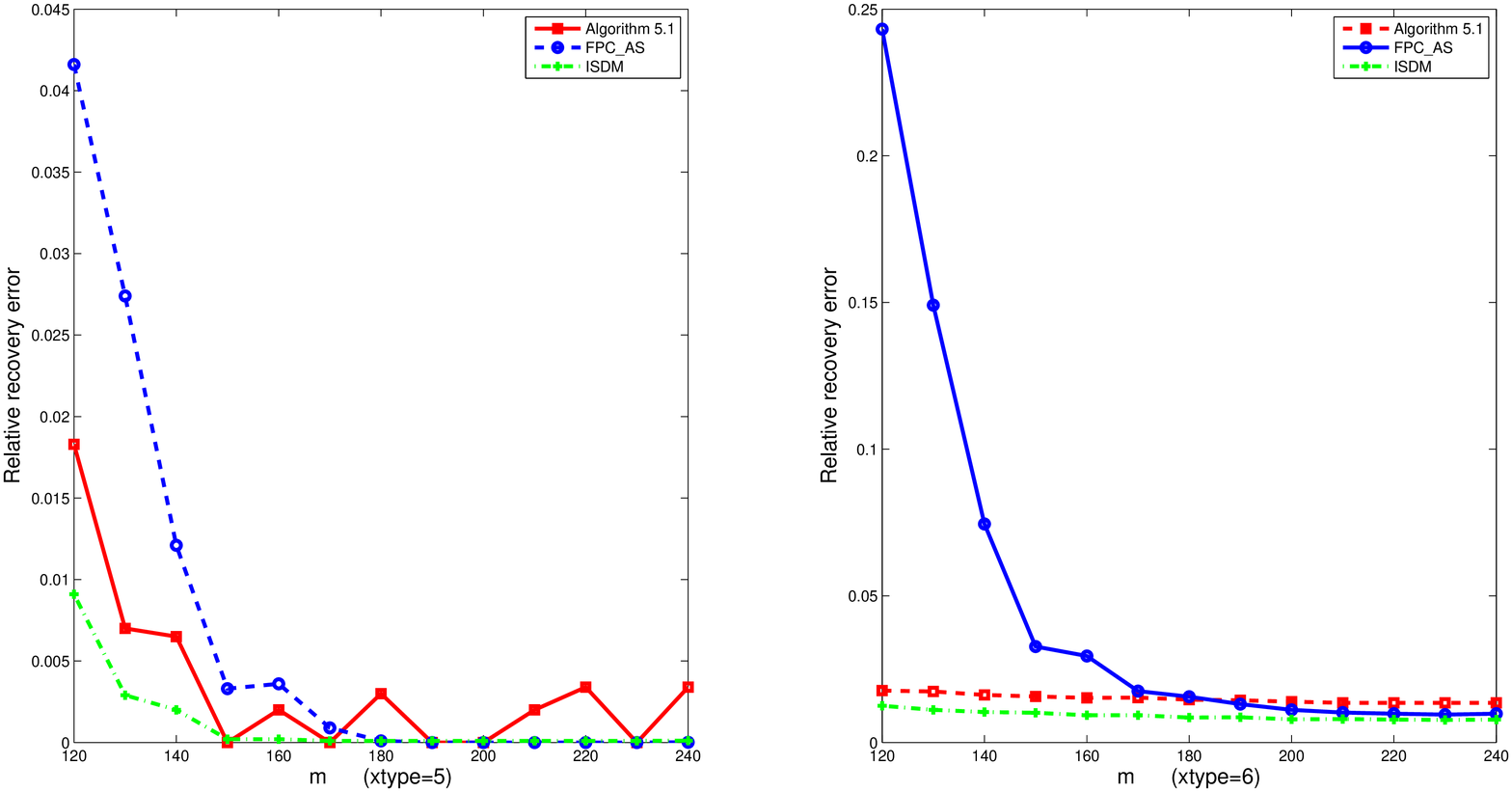}}
 \caption{\small Relative recovery error of three solvers (Atype$=3$)}
 \label{fig4}
 \end{center}
 \end{figure}

 \begin{figure}[htbp]
 \begin{center}
 {\includegraphics[width=0.95\textwidth]{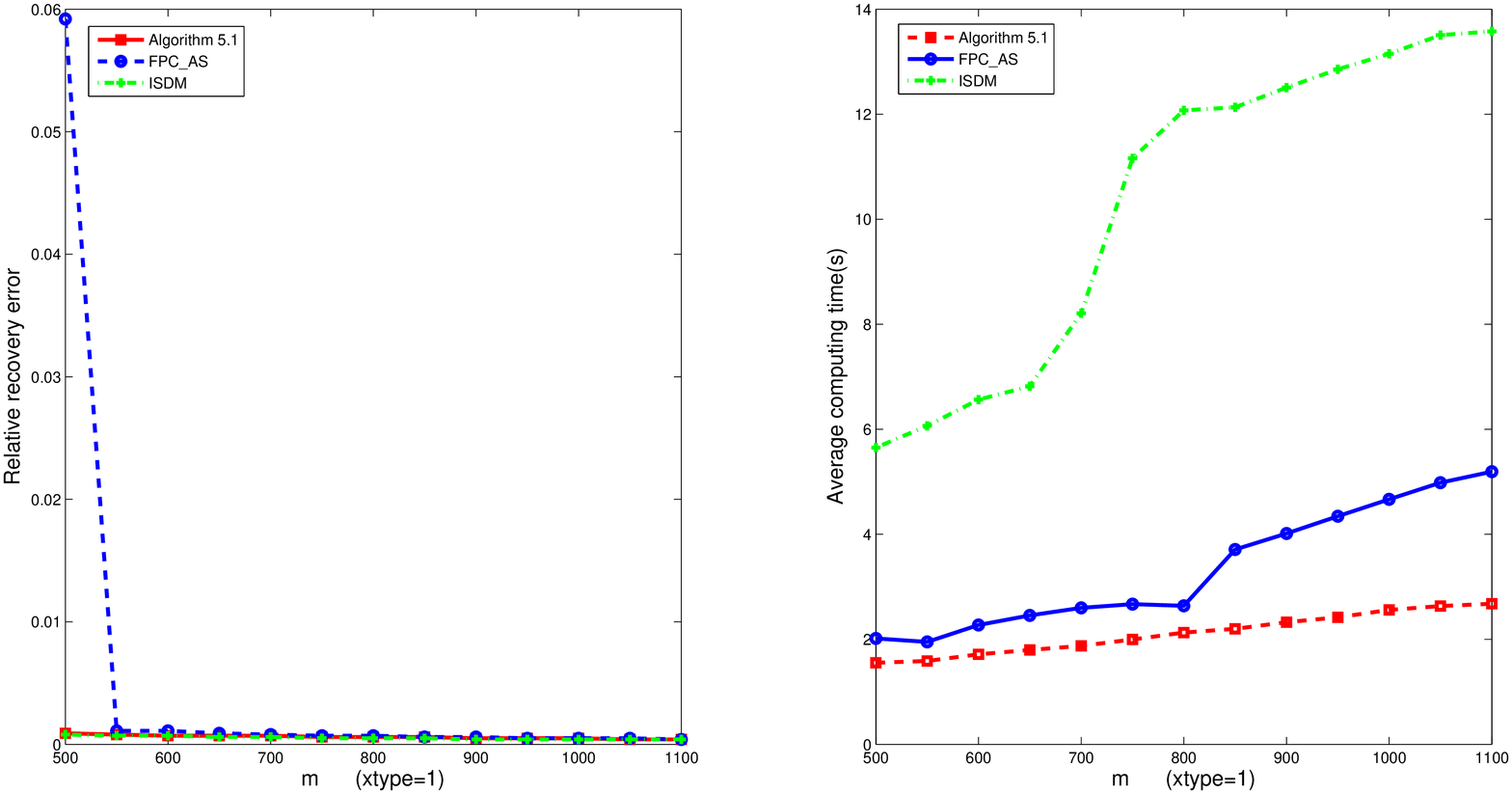}}\\
 {\includegraphics[width=0.95\textwidth]{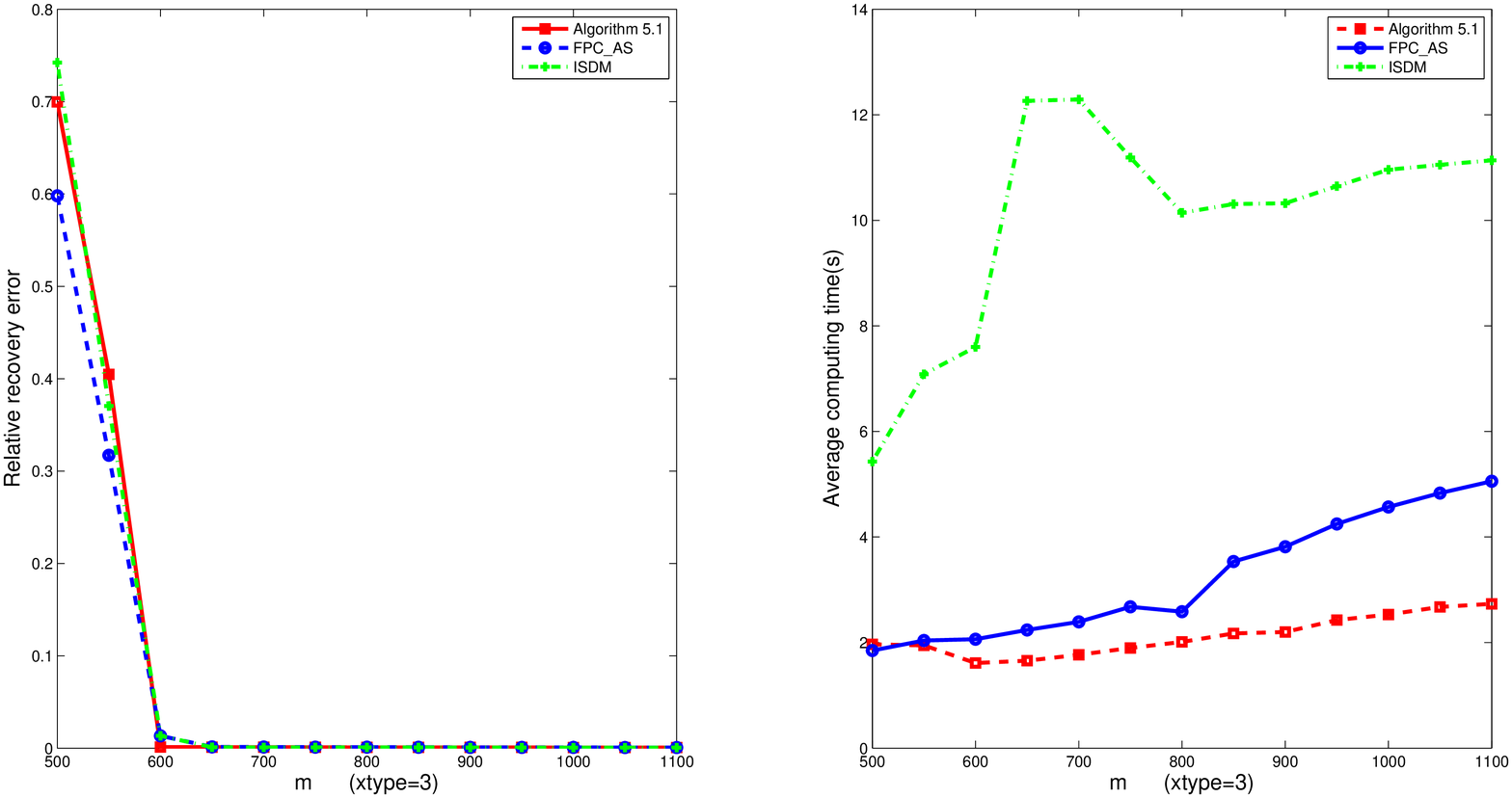}}\\
 {\includegraphics[width=0.95\textwidth]{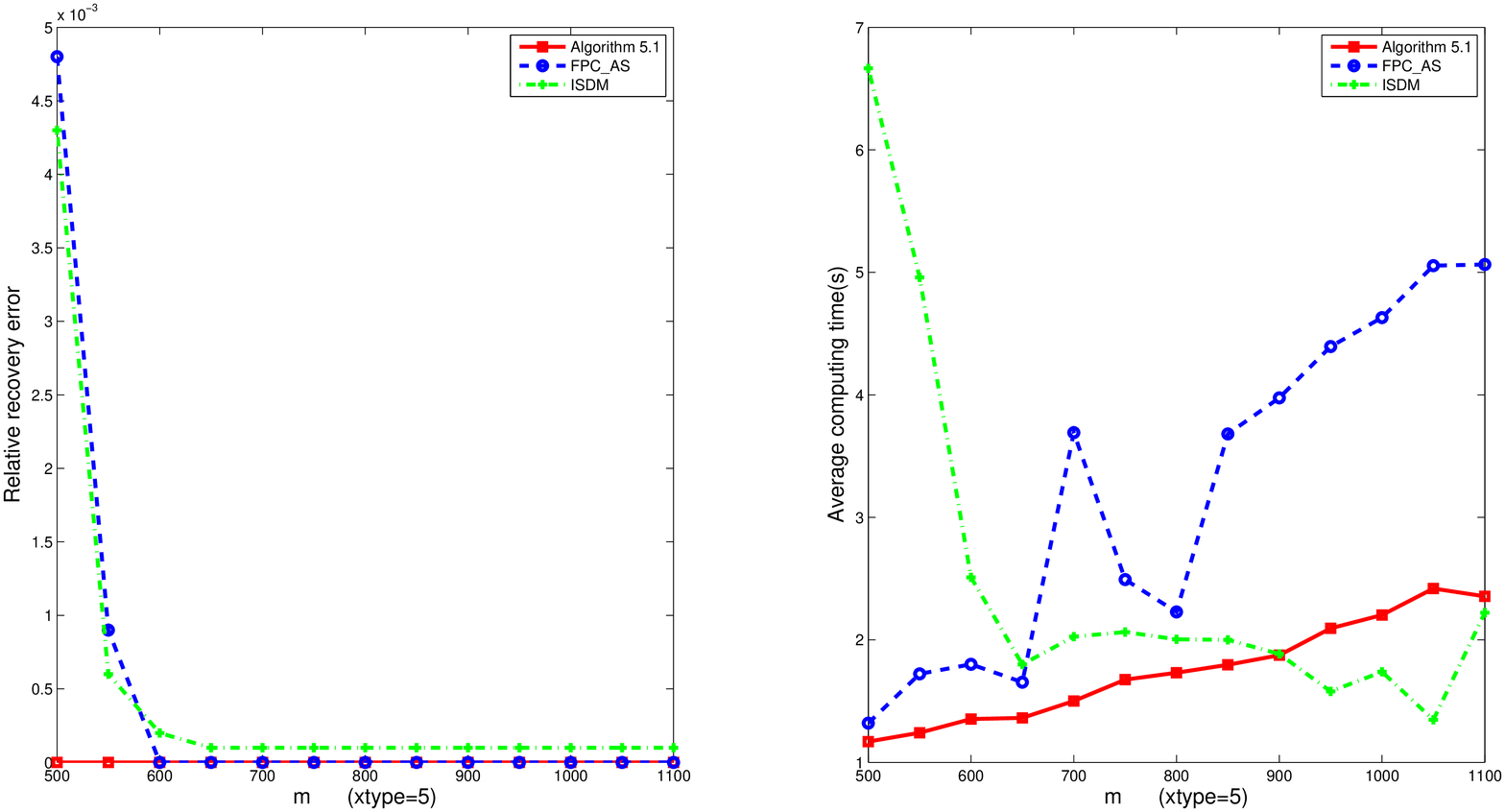}}
 \caption{\small Relative recovery error and average computing time of three solvers (Atype$=4$)}
 \label{fig5}
 \end{center}
 \end{figure}

  We first took the matrix of type $3$ for example to test the influence of the number of measurements
  $m$ on the recovery errors of three solvers for different types of signals.
  For each type of signals, we considered $n=600$ and $K=40$ and took the number of measurements
  $m\in\{120,130,\cdots,240\}$. For each $m$, we generated $50$ problems
  randomly and tested the recovery error of each solver. The curves in Figure \ref{fig4} depict
  how the relative recovery error of three solvers vary with the number of measurements for
  different types of signals. From this figure, we see that for the signals of types $1$, $2$ and $6$,
  Algorithm \ref{Hybrid-Alg1} and {\bf ISDM} require less measurements to yield the desirable
  recovery error than {\bf FPC\_AS} does; for the signals of types $3$ and $4$, the three solvers are comparable
  in terms of recovery errors; and for the signals of type $5$, {\bf ISDM} yields a little better recovery error
  than Algorithm \ref{Hybrid-Alg1} and {\bf FPC\_AS}. After checking, we found that for the signals of type $5$,
  the solutions yielded by {\bf ISDM} have a large average residual; for example, when $m=240$,
  the average residual attains $0.3$, whereas the average residual yielded by Algorithm \ref{Hybrid-Alg1}
  and {\bf FPC\_AS} are less than $0.02$. In other words, the solutions yielded by {\bf ISDM}
  deviate much from the set $\{x\in\Re^n\ |\ \|Ax-b\|\le \delta\}$.

  \medskip

  Finally, we took the matrix of type $4$ for example to compare the recovery errors and
  the computing time of three solvers for the signals of a little higher dimension.
  Since Figure \ref{fig4} shows that the three solvers have the similar performance for
  the signals of types $1$, $2$ and $6$, and the similar performance for the signals of type $3$ and $4$,
  we compared the performance of three solvers only for the signals of types $1$, $3$ and $5$.
  For each type of signals, we considered the dimension $n=2^{11}$ and the number of nonzeros $K=150$
  and took the number of measurements $m\in\{500,550,\cdots,1100\}$. For each $m$,
  we generated $50$ problems randomly and tested the recovery error of each solver.
  The curves of Figure \ref{fig5} depict how the recovery error and the computing time of
  three solvers vary with the number of measurements for different types of signals.

  \medskip

  Figure \ref{fig5} shows that for the signals of a little larger dimension,
  Algorithm \ref{Hybrid-Alg1} yields comparable recovery errors with {\bf ISDM} and {\bf FPC\_AS}
  and requires less computing time than {\bf ISDM} and {\bf FPC\_AS}. Together with Figure \ref{fig4},
  we conclude that for the noisy signal recovery, Algorithm \ref{Hybrid-Alg1} is superior to {\bf ISDM}
  and {\bf FPC\_AS} in terms of computing time, and comparable with {\bf ISDM} and better than {\bf FPC\_AS}
  in terms of the recovery error.

  \subsection{Sparco collection}\label{Subsec5.4}

  In this subsection we compare the performance of Algorithm \ref{Hybrid-Alg1} with that of
  {\bf FPC\_AS} and ${\bf ISDM}$ on ${\bf 24}$ problems from the Sparco collection \cite{BFHHS09},
  for which the matrix $A$ is stored implicitly. Table \ref{tab3} reports
  their numerical results where, each column has the same meaning as in Table \ref{tab2}.
  When the true $x^*$ is unknown, we mark {\bf Relerr} as ``--".

  \medskip

  From Table \ref{tab3}, we see that Algorithm \ref{Hybrid-Alg1} can solve those large-scale problems
  such as ``srcsep1", ``srcsep2", ``srcsep3", ``angiogram" and ``phantom2" with the desired feasibility,
  where ``srcsep3" has the dimension $n=196608$, and requires comparable computing time with
  {\bf FPC\_AS}, which is less than that required by {\bf ISDM} for almost all the test problems.
  The solutions yielded by Algorithm \ref{Hybrid-Alg1} have the smallest zero-norm for almost all test problems,
  and have better feasibility than those given by {\bf ISDM}. In particular, for those problems on
  which {\bf FPC\_AS} and {\bf ISDM} fail (for example, ``heavisgn'', ``blknheavi'' and ``yinyang''),
  Algorithm \ref{Hybrid-Alg1} still yields the desirable results. Also, we find that for some problems
  (for example, ``angiogram'' and ``phantom2''), the solutions yielded by {\bf FPC\_AS} have good feasibility,
  but their zero-norms are much larger than those of the solutions yielded by Algorithm \ref{Hybrid-Alg1} and {\bf ISDM}.

 \begin{table}[htbp]
 \begin{center}
 \caption{\label{tab3} Numerical comparisons of three solvers on Sparco collection}
 \begin{tabular}{|c|c|c|c|c|c|c|c|}
  \hline
   \ {\bf No.} &{\bf Problem} & {\bf Solver} & {\bf time(s)} & {\bf Relerr} & {\bf Res} & {\bf nMat}& {\bf nnzx}\\
  \hline
   \ &    & {\bf Algorithm} \ref{Hybrid-Alg1} & 0.45 & 2.23e-12 & 2.23e-10 & 618 & 4 \\
   \ 1 & {\rm Heavisine} & {\bf FPC\_AS} & 5.82 & 6.86e-1 & 1.54e+0 & 6891 & 1557\\
   \  &  & {\bf ISDM} & 5.55 &6.47e-1 & 2.17e-0 & -- & 10\\
    \hline
    \ &   & {\bf Algorithm} \ref{Hybrid-Alg1} & 0.14 & 2.53e-12 & 2.00e-10 & 539 & 71\\
    \ 2 & {\rm blocksig} & {\bf FPC\_AS} & 0.03 & 1.07e-11 & 8.43e-10 & 9 & 71\\
    \  &  & {\bf ISDM} & 0.03 & 7.94e-15 & 5.91e-13 &--  & 71\\
    \hline
    \ &   & {\bf Algorithm} \ref{Hybrid-Alg1} & 0.66 & 9.80e-13 & 9.95e-11 & 1546 & 115 \\
    \ 3& {\rm cosspike} & {\bf FPC\_AS} & 0.13 & 1.10e-11 & 1.10e-9 & 221 & 115\\
    \  &   & {\bf ISDM} & 5.99 & 4.25e-7 & 4.32e-5  &--  & 115\\
    \hline
     \  &   & {\bf Algorithm} \ref{Hybrid-Alg1} & 0.89 & 1.05e-10 & 1.09e-8 & 1748 & 121 \\
    \ 4& {\rm zsinspike} & {\bf FPC\_AS} & 4.74 & 4.32e-11 & 4.82e-9 & 6669 & 121\\
    \  &  & {\bf ISDM} & 5.13 & 3.68e-11 & 4.18e-9 & -- & 121\\
    \hline
    \  &  & {\bf Algorithm} \ref{Hybrid-Alg1} & 1.31 & 1.06e-8 & 6.59e-7 & 1730 & 59\\
    \ 5 & {\rm gcosspike} & {\bf FPC\_AS} & 2.07 & 1.48e-11 & 8.42e-10 & 1123 & 59\\
    \  &  & {\bf ISDM} & 8.81 & 4.34e-4 & 1.30e-2 & -- & 61 \\
    \hline
    \ &   & {\bf Algorithm} \ref{Hybrid-Alg1} & 168.0 & 5.63e-9 & 2.06e-5 & 3851  & 166\\
    \ 6 & {\rm p3poly} & {\bf FPC\_AS} & 23.4 & 5.85e-12 & 5.44e-11 & 1691 & 166 \\
    \  &  & {\bf ISDM} & 277.0 & 4.63e-3 & 9.05e-0 &--  & 211  \\
    \hline
    \ &   & {\bf Algorithm} \ref{Hybrid-Alg1} & 0.34 & 1.97e-11 & 4.33e-11 & 266 & 20 \\
    \ 7& {\rm sgnspike} & {\bf FPC\_AS} & 0.17 & 4.50e-10 & 9.36e-10 & 63  & 20\\
    \  &  & {\bf ISDM} & 0.81 & 1.95e-14 & 3.59e-14 & --  & 20 \\
    \hline
     \ &   & {\bf Algorithm} \ref{Hybrid-Alg1} & 1.37 & 5.04e-12 & 1.58e-11 & 609 & 20 \\
    \ 8 & {\rm zsgnspike} & {\bf FPC\_AS} & 33.3 & 3.57e-9 & 1.60e-8 & 6747 & 20\\
    \  &  & {\bf ISDM} & 3.15 & 1.90e-14 & 5.02e-14 &--  & 20  \\
    \hline
    \  &  & {\bf Algorithm} \ref{Hybrid-Alg1} & 0.19 & 3.03e-8 & 8.92e-7 & 2137  & 12\\
    \ 9& {\rm blkheavi} & {\bf FPC\_AS} & 0.14 & 2.45e-11 & 3.91e-10 & 789 & 12\\
    \  &  & {\bf ISDM} & 5.69 & 7.30e+2 & 1.02e+4 & --  & 102 \\
    \hline
    \  &  & {\bf Algorithm} \ref{Hybrid-Alg1} & 0.45 & 3.29e-7 & 2.00e-6 & 2059 & 12\\
    \ 10 & {\rm blknheavi} & {\bf FPC\_AS} & 1.72 & 2.57e-2 & 1.40e-1 & 6797 & 344\\
    \   & & {\bf ISDM} & 2.17 & 6.05e-1 & 1.54e+0  &--  & 25\\
    \hline
    \  &  & {\bf Algorithm} \ref{Hybrid-Alg1} & 0.31 & 1.64e-9 & 1.46e-7 & 1671  & 32 \\
    \ 11 & {\rm gausspike} & {\bf FPC\_AS} & 0.14 & 3.16e-12 & 4.38e-11 & 181 & 32\\
    \  &  & {\bf ISDM} & 0.99 & 2.89e-14 & 1.67e-12  & -- & 32 \\
    \hline
     \ &    & {\bf Algorithm} \ref{Hybrid-Alg1} & 186.0 & -- & 6.81e-6 & 3431 & 21520 \\
    \ 12& {\rm srcsep1} & {\bf FPC\_AS} & 238.0 & -- & 4.08e-5  & 6885 & 42676\\
    \  &  & {\bf ISDM} & 388.0 & -- & 3.81e-3  & -- & 21644 \\
    \hline
     \ &   & {\bf Algorithm} \ref{Hybrid-Alg1} & 380.0 & -- & 1.92e-6 & 3628 & 21733\\
    \ 13& {\rm srcsep2} & {\bf FPC\_AS} & 351.0 & -- & 3.09e-4  & 6885 & 64478\\
    \  &  & {\bf ISDM} & 601.0 & -- & 2.01e-3 & -- & 23258 \\
    \hline

  \end{tabular}
  \end{center}
  \end{table}
  \begin{table}[htbp]
  \begin{center}
  \begin{tabular}{|c|c|c|c|c|c|c|c|}
   \hline
    \  &  & {\bf Algorithm} \ref{Hybrid-Alg1} & 409.0 & -- & 1.83e-4 & 3934 & 110406  \\
    \ 14& {\rm srcsep3} & {\bf FPC\_AS} & 589.0 & -- & 1.25e-7 & 7131 & 113438 \\
    \   & & {\bf ISDM} & 316.0 & -- & 3.78e-3  & -- & 110599 \\
    \hline
    \  &  & {\bf Algorithm} \ref{Hybrid-Alg1} & 3.53 & -- & 3.31e-6 & 2399 & 606\\
    \ 15 & {\rm phantom1} & {\bf FPC\_AS} & 19.4 & -- & 1.19e-5  & 7203 & 3989\\
    \  &  & {\bf ISDM} & 9.08 & -- & 2.22e-10 & -- & 811  \\
    \hline
    \  &  & {\bf Algorithm} \ref{Hybrid-Alg1} & 1.59 & -- & 6.96e-7 & 695 & 574\\
    \ 16 & {\rm angiogram} & {\bf FPC\_AS} & 4.24 & -- & 3.29e-9 & 485 & 9881\\
    \  &  & {\bf ISDM} & 3.34 & -- & 2.17e-13  & -- & 574\\
    \hline
    \  &  & {\bf Algorithm} \ref{Hybrid-Alg1} & 57.8 & -- & 2.97e-7 & 2914 & 20962\\
    \ 17& {\rm phantom2} & {\bf FPC\_AS} & 44.4 & -- & 1.04e-8 & 861 & 64599\\
    \  &  & {\bf ISDM} & 149.0 & -- & 5.93e-7 & -- & 29313  \\
     \hline
    \  &  & {\bf Algorithm} \ref{Hybrid-Alg1} & 53.0 & -- & 6.41e+1  & 5102 & 1831\\
    \ 18& {\rm smooth soccer} & {\bf FPC\_AS} & 97.8 & -- & 5.18e+0 & 6875  & 3338\\
    \  &  & {\bf ISDM} & 456.0 & -- & 2.64e+3 & --  & 2\\
        \hline
    \  &  & {\bf Algorithm} \ref{Hybrid-Alg1} & 35.8 & -- & 1.45e-6 & 4427 & 701\\
    \ 19& {\rm soccer} & {\bf FPC\_AS} & 25.7 & -- & 4.61e-10& 1769  & 701 \\
    \  &  & {\bf ISDM} & 247.0 & -- & 1.00e+7  & -- & 6\\
    \hline
    \  &  & {\bf Algorithm} \ref{Hybrid-Alg1} & 8.70 & -- & 5.65e-7 & 2710 & 771\\
    \ 20 & {\rm yinyang} & {\bf FPC\_AS} & 35.4 & -- & 7.99e-3 & 6563  & 3281\\
    \  &  & {\bf ISDM} & 52.7 & -- & 6.34e-4 & -- & 886 \\
    \hline
    \  &  & {\bf Algorithm} \ref{Hybrid-Alg1} & 167.0 & -- & 4.43e-6  & 3992 & 62757\\
    \ 21 & {\rm blurrycam} & {\bf FPC\_AS} & 76.5 & -- & 8.63e-7  & 2313 & 62757\\
    \  &  & {\bf ISDM} & 525.0 & -- & 2.51e-1 & --  & 54829 \\
    \hline
    \  &  & {\bf Algorithm} \ref{Hybrid-Alg1} & 27.6 & -- & 1.36e-6  & 3225 & 15592\\
    \ 22& {\rm blurspike} & {\bf FPC\_AS} & 11.7 & -- & 4.86e-7  & 1863 & 15592\\
    \   & & {\bf ISDM} & 67.0 & -- & 2.62e-3 & --  & 15276\\
    \hline
     \  &  & {\bf Algorithm} \ref{Hybrid-Alg1} & 0.03 & 6.06e-10 & 3.02e-10 & 50 & 3 \\
    \ 23 & {\rm jitter} & {\bf FPC\_AS} &0.02 & 7.99e-10 & 3.84e-10 & 35 & 3 \\
    \  &  & {\bf ISDM} & 0.02 & 2.10e-14 & 9.69e-15 & -- & 3 \\
    \hline
    \  &  & {\bf Algorithm} \ref{Hybrid-Alg1} & 0.36 & 3.48e-9 & 1.27e-7  & 1851 & 12\\
    \ 24 & {\rm spiketrn} & {\bf FPC\_AS} & 1.45 & 3.81e-11 & 1.94e-9  & 4535 & 12\\
    \  &  & {\bf ISDM} & 4.82 & 3.68e-0 & 1.40e+1  & -- & 34\\
    \hline
  \end{tabular}
  \end{center}
  \end{table}

   \medskip

   From the numerical comparisons in Subsection \ref{Subsec5.1}-\ref{Subsec5.4}, we conclude that
   Algorithm \ref{Hybrid-Alg1} is comparable even superior to {\bf ISDM} in terms of recoverability,
   and the superiority of Algorithm \ref{Hybrid-Alg1} is more remarkable for
   those difficult problems from Sparco collection. The recoverability of Algorithm \ref{Hybrid-Alg1} and {\bf ISDM} is
   higher than that of {\bf FPC\_AS}. In particular, Algorithm \ref{Hybrid-Alg1} requires
   less computing time than {\bf ISDM}. The recoverability
   and recovery error of {\bf QPDM} is much worse than that of the other three solvers.

  \section{Conclusions}

   In this work we reformulated the zero-norm problem (\ref{l0}) as
   an equivalent MPEC, then established its exact penalty formulation (\ref{MPEC1}).
   To the best of our knowledge, this novel result can not be obtained
   from the existing exact penalty results for MPECs. Motivated by the special
   structure of exact penalty problem, we proposed a decomposition method
   for dealing with the MPEC problem, and consequently the zero-norm problem.
   This method consists of finding the solution of a finite number of
   weighted $l_1$-norm minimization problems, for which we propose an effective partial PPA algorithm
   for dealing with them. In particular, we show that this method can yield an optimal solution of
   the zero-norm problem under the null space condition used in \cite{KXAB11}.
   Numerical comparisons show that the exact penalty decomposition method  is significantly
   better than the quadratic penalty decomposition method \cite{ZL10}, is comparable with
   {\bf ISDM} in terms  of recoverability \cite{WY10} but requires less computing time,
   and has better recoverability than {\bf FPC\_AS} \cite{WYGZ10} and requires comparable computing time.

  \medskip

  There are several research topics worthwhile to pursue; for example, one may consider to
  extend the results of this paper to rank minimization problems, design other effective
  convex relaxation methods for (\ref{l0}) based on its equivalent MPEC problem, and make
  numerical comparisons for the exact penalty decomposition method with the weighted $l_1$-norm
  subproblems solved by different effective algorithms.

 \bigskip
 \noindent
 {\bf Acknowledgements.} We would like to thank Professor Jiye Han
  for giving us some helpful comments on the closeness of the image
  of a closed convex cone under a linear mapping, Professor Defeng Sun
  at  National University of Singapore for some important suggestions
  on the solution of subproblems and the revision of this paper,
  and Dr. Weimin Miao for helpful discussions on the null space properties
  of a linear operator. Last but not least, thanks also go to the anonymous referees
  for carefully reading our paper and providing valuable suggestions for its revision.


 \end{document}